\documentclass{elsarticle}
\usepackage{mystyle, multirow}
\usepackage[margin=2.5cm]{geometry}

\makeatletter                   
\def\mdseries@tt{m}             
\makeatother                    

\usepackage[draft=true]{minted}
\usepackage{xcolor,cancel}

\definecolor{lightgray}{gray}{.95}
\newminted{python}{
    linenos=true,
    bgcolor=lightgray,
    tabsize=4,
    fontfamily=courier,
    fontsize=\small,
    xleftmargin=5pt,
    xrightmargin=5pt
}

\DeclareMathOperator*{\cotan}{cotan}
\DeclareMathOperator*{\dz}{dz}
\DeclareMathOperator*{\dev}{dev}
\DeclareMathOperator*{\sym}{sym}
\renewcommand{\skew}{\operatorname{skew}}

\graphicspath{{./pic/},{./results/}}
\renewcommand{\udl}[1]{\boldsymbol{#1}}
\renewcommand{\uddl}[1]{\boldsymbol{#1}}

\newcommand*{\mydots}{\vcenter{\offinterlineskip \vskip-.1ex\hbox{$.$}\vskip-.7ex\hbox{$.$}\vskip-.7ex\hbox{$.$}}} 

\newcommand{\tripleC}{\mathbin{{\mathchoice{\textstyle\mydots}
{\textstyle\mydots} {\scriptstyle\scalebox{0.8}{$\mydots$}} {\scriptscriptstyle\scalebox{0.6}{$\mydots$}}}}}
\begin{document}
\sloppy

\begin{frontmatter}
\title{Automated formulation and resolution of limit analysis problems}
\author[myu]{Jeremy Bleyer \corref{cor1}}
\ead{jeremy.bleyer@enpc.fr}
\author[myu]{Ghazi Hassen}
\ead{ghazi.hassen@enpc.fr}
\address[myu]{Laboratoire Navier, Ecole des Ponts ParisTech, Univ Gustave Eiffel, CNRS \\
  6-8 av Blaise Pascal, Cit\'e Descartes, 77455 Champs-sur-Marne, FRANCE}
\cortext[cor1]{Corresponding author}

\begin{abstract}
The present manuscript presents a framework for automating the formulation and resolution of limit analysis problems in a very general manner. This framework relies on FEniCS domain-specific language and the representation of material strength criteria and their corresponding support function in the conic programming setting. Various choices of finite element discretization, including discontinuous Galerkin interpolations, are offered by FEniCS, enabling to formulate lower bound equilibrium elements or upper bound elements including discontinuities for instance. The numerical resolution of the corresponding optimization problem is carried out by the interior-point solver \texttt{Mosek} which takes advantage of the conic representation for yield criteria. Through various illustrative examples ranging from classical continuum limit analysis problems to generalized mechanical models such as plates, shells, strain gradient or Cosserat continua, we show that limit analysis problems can be formulated using only a few lines of code, discretized in a very simple manner and solved extremely efficiently. This paper is accompanied by a FEniCS toolbox implementing the above-mentioned framework.
\end{abstract}

\end{frontmatter}

\section{Introduction}

Limit analysis \cite{hill1950mathematical}, or more generally yield design theory \cite{salencon1983calcul,salenccon2013yield}, is an efficient method for computing the ultimate load, or bearing capacity, of a structure based on the sole knowledge of a given local strength criterion (or plasticity criterion for limit analysis) and applied external loads. Its main advantage is that it is a direct method i.e. it avoids the need of performing a step-by-step non-linear analysis until complete failure, which is estimated poorly in general since it is associated with non-convergence of the resolution algorithm. More precisely, limit analysis theory can be formulated as a convex optimization problem and therefore benefits from variational approximations on the primal and associated dual problem yielding the so-called lower bound static and upper bound kinematic approaches. The exact collapse load can therefore be bracketed by the bounding status of the static and kinematic solutions. This method has found tremendous applications in mechanical and civil engineering problems since analytical upper bounds can be obtained very efficiently by considering simple collapse mechanisms but also because lower bounds provide a safe approximation to the exact collapse load. Typical fields of application include soil slope stability, footing bearing capacity or other geotechnical problems \cite{chen2013limit}, rigid-block masonry structures, design of reinforced-concrete structures \cite{chen2007plasticity, nielsen2016limit}, especially through strut-and-tie methods \cite{schlaich1987toward} or yield-line analysis \cite{johansen1962yield}, collapse loads of frame, plate or shell structures \cite{save1995atlas,save1997plastic}, etc.

In all the above-mentioned applications, hand-based solutions are quite easy to compute and can be found in different design codes. However, limit analysis techniques have been somehow limited to hand-based solutions for quite a long time because of the difficulties encountered in the past when automating their resolution in a finite-element discrete setting for instance. Indeed, the corresponding optimization problems are inherently non-smooth and large-scale which makes them challenging to solve. Standard gradient-based optimizers for instance are not a good candidate for their resolution due to the highly non-smooth aspect. The major breakthrough in this field is associated with the development of efficient optimization algorithms particularly well-suited for this class of optimization problems. The first progress has been made with the simplex algorithm \cite{dantzig1955generalized} for solving linear programming (LP) problems. However, limit analysis problems fall into the LP category only when the strength criterion is a polytope i.e. when it can be written as a collection of linear inequalities \cite{anderheggen1972finite,sloan1988lower}, whereas the majority of strength/plasticity criteria are non-linear. The second progress has been achieved with the advent of so-called \textit{interior-point} (IP) solvers which improved the complexity of solving LP problems \cite{karmarkar1984new}, the simplex algorithm has exponential complexity compared to polynomial (and, in practice, quasi-linear) complexity for interior-point solvers. More importantly, the interior-point method has been extended to more complex convex programs \cite{nesterov1994interior} such as second-order cone programs (SOCP) or semi-definite programs (SDP). Combining the weak complexity of IP solvers and the fact that most strength criteria can be expressed using second-order cone or semi-definite constraints \cite{bisbos2007second,makrodimopoulos2010remarks} now enables to solve complex and large-scale limit analysis problems \cite{krabbenhoft2008three,martin_finite-element_2008,portioli2014limit,vincent2018yield}. Obviously, a few other alternative methods have been proposed but, to the authors' opinion, none of them have been shown to exhibit a performance similar to state-of-the-art IP solvers.

As regards numerical discretization techniques for limit analysis problem, the vast majority of works relies on the finite-element method \cite{lysmer1970limit,anderheggen1972finite,pastor1976mise}. The specificity of limit analysis problems, compared to more standard nonlinear computations with displacement-based FE discretizations, lies in the use of static equilibrium-based finite-elements for obtaining true lower bound  \cite{sloan1988lower,krabbenhoft2002lower} (and therefore safe) estimates of the limit load but also in the use of discontinuous finite-elements for the kinematic upper bound approach \cite{sloan-kleeman1994discontinuous,makrodimopoulos2008upper,makrodimopoulos2007upper}. Indeed, most hand-based upper bound solutions have been obtained considering rigid-block mechanisms, thus involving no deformation but only displacement jumps in the plastic dissipation computation. Despite the higher computational cost compared to equivalent continuous interpolations, discontinuous interpolations provide more accurate limit load estimates \cite{pastor1976mise,sloan-kleeman1994discontinuous,krabbenhoft2005new,makrodimopoulos2008upper}, especially if finite-element edges are well-oriented. They are also more robust for certain problems since they do not suffer from locking issues, see for instance \cite{nagtegaal1974numerically} for volumetric locking in pressure-insensitive materials or \cite{bleyer2016numerical} for shear-locking in thin plates.

Due to the specific nature of the optimization problems, formulating a discretized version of a limit analysis (either static or kinematic) approach requires forming matrices representing, for instance, equilibrium, continuity or boundary conditions but also other linear relations coupling mechanical variables (such as stress or strain) with auxiliary variables used to express the strength conditions in a LP/SOCP or SDP format. Besides, depending on the specific choice made for the optimization solver, the standard input format of the problem may differ. As a result, discrete limit analysis problems require access to matrices that are not readily available from standard displacement-based FE solver and must therefore be implemented in an external program before calling the optimizer. Combining this aspect with the various types of FE discretizations and mechanical models makes the automation of limit analysis problem a challenge. As a result, limit analysis codes are usually limited to specific situations, sometimes with specific strength criteria.\\

In the present manuscript, we describe a general framework for the formulation of limit analysis problems for different mechanical models (2D/3D continua, plates/shells or generalized continua). Relying on the FEniCS finite-element library and symbolic representation of operators and code-generation capabilities, different FE discretization schemes (including discontinuous or equilibrium elements) can be easily formulated and generalized to more advanced mechanical models. The proposed framework therefore offers four levels of generality in the problem formulation:
\begin{itemize}
\item choice of a \textit{mechanical model}: limit analysis problems possess the same structure and can be formulated in a symbolic fashion through generalized stress/strain definitions (section \ref{sec:generalized-LA}).
\item choice of a \textit{strength criterion}: the formulation of its conic representation at a local point suffices to completely characterize the strength criterion, its translation to the global optimization problem being automatically performed (section \ref{sec:conic-repr}).
\item choice of a \textit{FE discretization}: including element type, interpolation degree or quadrature rule, all compatible for a variable number of degrees of freedom related to the choice of the mechanical model (section \ref{sec:continuum-LA}).
\item choice of the \textit{optimization solver}: although the accompanying FEniCS toolbox relies extensively on the \texttt{Mosek} optimization solver \cite{mosek}, once formulated in a standard conic programming form, the problem can then be written in a specific file format appropriate for another solver.
\end{itemize}
This versatility is further illustrated by considering homogenization theory in a limit analysis setting (section \ref{sec:homog}), plate and shell problems (section \ref{sec:plates}) and generalized continua such as Cosserat or strain gradient models (section \ref{sec:gen-cont}).

As regards numerical implementation, the present paper is accompanied by a Python module \texttt{fenics\_optim.limit\_analysis} implemented as a submodule of the FEniCS convex optimization package \texttt{fenics\_optim} \cite{jeremy_bleyer_2020_3604086} described in \cite{bleyer2019automating}.

\paragraph{Notations}
$\boldsymbol{A}:\boldsymbol{B}=A_{ij}B_{ji}$, $\boldsymbol{\Aa}\tripleC\boldsymbol{\Bb}=\Aa_{ijk}\Bb_{kji}$.

\section{\label{sec:generalized-LA}A general framework for limit analysis problems} 
In this section, we consider a material domain $\Omega \subset \RR^{d}$ (with $d=1,2,3$) associated with a specific mechanical model. In the subsequent applications, we will consider classical continuum theories such as 2D or 3D Cauchy continua or Reissner-Mindlin plate models for instance but also generalized continuum models encompassing higher-grade or higher-order theories. For this reason, the subsequent presentation will be made in a generalized continuum framework in which the mechanical stress or strain measures, equilibrium or continuity equations and boundary conditions will be written in an abstract fashion, their precise expression remaining to be specified for each particular mechanical theory. In particular, the presentation will be in the line of Germain's construction through the virtual work principle \cite{germain1973method,germain1973methode}.

\subsection{Virtual work principle for generalized continua}
Let us therefore consider a generalized virtual velocity field $\bu(\bx)$ of dimension $n$ and a set of strain measures $\bD \bu$ of dimension $m$ with $\bD$ being a generalized strain operator. Following \cite{germain1973method,germain1973methode}, such strain measures must be objective i.e. null for any rigid body motion. The virtual power of internal forces is assumed to be given by an internal force density depending linearly upon the strain measures:
\begin{equation}
\Pp^{(i)}(\bu)= -\int_\Omega \bSig\cdot \bD\bu \dx
\end{equation}
in which $\bSig$ denotes the generalized stress measure associated with $\bD\bu$ by duality. The above expression must in fact be understood in the sense of distributions i.e. $\bu$ may exhibit discontinuities $\bJ \bu$ (consistent with the definition of operator $\bD$) across some internal surface $\Gamma$. The power of internal forces therefore writes more explicitly as:
\begin{equation}
\Pp^{(i)}(\bu)= -\left(\int_{\Omega\setminus\Gamma} \bSig\cdot \bD\bu \dx + \int_\Gamma \bSig\cdot\bJ \bu \dS\right) 
\end{equation}
The power of external forces is assumed to consist of two contributions: long-range interaction forces described by a volume density $\boldsymbol{F}$ and boundary contact forces described by a surface density $\boldsymbol{T}$ acting on the exterior boundary $\partial \Omega$. Each power density depends linearly upon the generalized velocity so that the total power is given by:
\begin{equation}
\Pp^{(e)}(\bu)= \int_{\Omega} \boldsymbol{F}\cdot \bu \dx + \int_{\partial \Omega} \boldsymbol{T}\cdot \bu \dS 
\end{equation}

According to the virtual power principle, the system is in equilibrium if and only if the sum of the internal and external virtual powers is zero for any virtual velocity field:
\begin{equation}
\Pp^{(i)}(\bu)+\Pp^{(e)}(\bu) = 0 \quad \forall \bu \tag{\texttt{equilibrium}}
\label{eq:equilibrium}
\end{equation}

\subsection{General formulation of a limit analysis problem}

Limit analysis (or yield design) theory amounts to finding the maximum loading a system can sustain considering only equilibrium and strength conditions for its constitutive material. The latter can be generally described by the fact that the generalized stresses $\bSig(\bx)$ must belong to a strength domain $G(\bx)$ for all point $\bx\in\Omega$:
\begin{equation}
\bSig(\bx)\in G(\bx) \quad \forall \bx\in\Omega \hspace{-1.1cm}\tag{\texttt{strength condition}}
\label{eq:crit}
\end{equation}
The strength domain $G \subseteq \RR^m$ is assumed to be a convex set (it may be unbounded and non-smooth) which usually contains the origin $0\in G$.\\

Finding the maximum loading will be achieved with respect to a given loading \textit{direction} i.e. by assuming that both volume and surface forces
depend upon a single load factor $\lambda$ in an affine manner:
\begin{align}
\boldsymbol{F}(\lambda) &= \lambda\boldsymbol{f}+\boldsymbol{f}_0\\
\boldsymbol{T}(\lambda) &= \lambda\boldsymbol{t}+\boldsymbol{t}_0
\end{align}
where $\boldsymbol{f},\boldsymbol{t}$ are given loading directions and $\boldsymbol{f}_0,\boldsymbol{t}_0$ are fixed reference loads. In such a case, the power of external forces can also be written in an affine manner with respect to $\lambda$:
\begin{equation}
\Pp^{(e)}_\lambda(\bu) = \lambda\Pp^{(e)}(\bu) + \Pp^{(e)}_0(\bu) 
\end{equation}
with obvious notations. Let us mention that if one wants to describe the set of ultimate load defined by multiple loading parameters (see for instance section \ref{sec:homog}), the loading directions must be varied in order to describe the boundary of the ultimate load domain.\\

The limit analysis problem can finally be formulated as finding the maximum load factor $\lambda$ such that there exists a generalized stress field $\bSig(\bx)$ in equilibrium with $(\boldsymbol{F}(\lambda),\boldsymbol{T}(\lambda))$ and complying with the material strength properties i.e. satisfying both \eqref{eq:equilibrium} and \eqref{eq:crit} which can also be written as:
\begin{equation}
\begin{array}{rl}
\displaystyle{\lambda^+ = \sup_{\lambda,\bSig}} & \lambda \\
\text{s.t.} & \Pp^{(i)}(\bu)+\lambda\Pp^{(e)}(\bu)+\Pp^{(e)}_0(\bu) = 0 \quad \forall \bu \\
& \bSig(\bx) \in G(\bx) \quad \forall \bx \in \Omega
\end{array}
\label{car-mixed}
\end{equation}
Let us mention that for the infinite-dimensional convex problem \eqref{car-mixed} to have a solution, the fixed loading $(\boldsymbol{f}_0,\boldsymbol{t}_0)$ must be a sustainable loading i.e. there must exist a stress field in equilibrium with $(\boldsymbol{f}_0,\boldsymbol{t}_0)$ and satisfying \eqref{eq:crit}.

Formulation \eqref{car-mixed} will be the basis of the mixed finite-element formulation discussed in section \ref{sec:mixed-FE} when choosing proper interpolation spaces for $\bSig$ and $\bu$. We now turn to the general formulation of the static and kinematic approaches.

\subsubsection{Static approach}
Starting from the weak formulation of equilibrium given by \eqref{eq:equilibrium}, strong balance equations, continuity conditions and boundary conditions can be obtained for the generalized stresses $\bSig$. These will generally take the following form:
\begin{align}
\boldsymbol{\mathcal{E}}\bSig + \lambda \boldsymbol{f}+\boldsymbol{f}_0 = 0 & \quad\text{in }\Omega \\
\boldsymbol{\mathcal{C}}\bSig = 0 & \quad\text{on }\Gamma\\
\boldsymbol{\mathcal{T}}\bSig = \lambda\boldsymbol{t}+\boldsymbol{t}_0 & \quad\text{on } \partial \Omega
\end{align}
where $\boldsymbol{\mathcal{E}}$ is an equilibrium operator (adjoint to $\boldsymbol{\mathcal{D}}$) and $\boldsymbol{\mathcal{C}}$ and $\boldsymbol{\mathcal{T}}$ are some continuity and trace operators related to $\boldsymbol{\mathcal{E}}$. A generalized stress field $\bSig(\bx)$ satisfying these conditions will be termed as \textit{statically admissible} with a given loading $(\lambda\boldsymbol{f}+\boldsymbol{f}_0, \lambda\boldsymbol{t}+\boldsymbol{t}_0)$.

The pure static formulation can therefore be generally written as:
\begin{equation}
\begin{array}{rl}
\displaystyle{\lambda^+ = \sup_{\lambda,\bSig}} & \lambda \\
\text{s.t.} & \boldsymbol{\mathcal{E}}\bSig + \lambda \boldsymbol{f}+\boldsymbol{f}_0 = 0 \quad\text{in }\Omega \\
& \boldsymbol{\mathcal{C}}\bSig = 0 \quad\text{on }\Gamma \\
& \boldsymbol{\mathcal{T}}\bSig = \lambda\boldsymbol{t}+\boldsymbol{t}_0 \quad\text{on } \partial \Omega\\
& \bSig(\bx) \in G(\bx) \quad \forall \bx \in \Omega
\end{array}\tag{\texttt{SA}}
\label{car-stat}
\end{equation}
Obviously $\bSig$ must belong to an appropriate functional space $\mathcal{W}$ consistent with the nature of the above operators. If one restricts to a (finite-dimensional) subset $\mathcal{W}_h \subset \mathcal{W}$ such that all constraints of \eqref{car-stat} can be satisfied exactly, the corresponding solution $\lambda_s$ of the corresponding (finite) convex optimization problem will therefore be a lower bound to the exact limit load: $\lambda_s \leq \lambda^+$.

\subsubsection{Kinematic approach}
The kinematic formulation of a limit analysis problem can be obtained from an equivalent formulation of \eqref{car-mixed}:
\begin{equation}
\lambda \leq \lambda^+ \Longleftrightarrow \exists \bSig \text{ s.t. } \begin{cases}
\lambda\Pp^{(e)}(\bu)+\Pp^{(e)}_0(\bu) = - \Pp^{(i)}(\bu) \quad \forall \bu  \\
\bSig(\bx)\in G(\bx) \quad \forall \bx\in \Omega
\end{cases}
\end{equation}
One therefore has that:
\begin{equation}
\lambda\Pp^{(e)}(\bu)+\Pp^{(e)}_0(\bu) \leq \sup_{\bSig(\bx) \in G(\bx)}\{-\Pp^{(i)}(\bu)\} = \Pp^{(mr)}(\bu) \quad \forall \bu
\end{equation}
where we introduced the maximum resisting power as:
\begin{equation}
\Pp^{(mr)}(\bu) = \int_{\Omega\setminus\Gamma} \pi_G(\bD\bu) \dx + \int_\Gamma \pi_G(\bJ \bu) \dS
\end{equation}
in which $\pi_G$ is the support function of the convex set $G$:
\begin{equation}
\pi_G(\boldsymbol{d}) = \sup_{\bSig\in G} \{\bSig\cdot \boldsymbol{d}\}
\end{equation}
One can conclude that:
\begin{equation}
\lambda \leq \lambda^+ \Longrightarrow \lambda \leq \dfrac{\Pp^{(mr)}(\bu)-\Pp^{(e)}_0(\bu)}{\Pp^{(e)}(\bu)} \quad \forall \bu
\end{equation}
Minimizing the right-hand side of the above relation therefore gives an upper bound to the exact maximal load $\lambda^+$. Under appropriate mathematical assumptions \cite{nayroles1970essai,fremond1978analyse} (non-restrictive in practice), it can be shown that the minimum is in fact $\lambda^+$ so that one has:
\begin{equation}
\lambda^+ = \inf_{\bu} \dfrac{\Pp^{(mr)}(\bu)-\Pp^{(e)}_0(\bu)}{\Pp^{(e)}(\bu)}
\end{equation} 
Observing that the above quotient is invariant when rescaling $\bu$ by a positive factor, a normalization constraint can in fact be considered to remove the denominator so that, one finally has for the kinematic approach:
\begin{equation}
\begin{array}{rl}
\displaystyle{\lambda^+ = \inf_{\bu}} & \Pp^{(mr)}(\bu)-\Pp^{(e)}_0(\bu) \\
\text{s.t.} & \Pp^{(e)}(\bu)=1
\end{array} \tag{\texttt{KA}}
\label{car-kin}
\end{equation}

Similarly to the static approach \eqref{car-stat}, $\bu$ must belong to an appropriate functional space $\mathcal{V}$. If one restricts to a (finite-dimensional) subset $\mathcal{V}_h \subset \mathcal{V}$, the corresponding solution $\lambda_u$ of the corresponding (finite) convex optimization problem will therefore be an upper bound to the exact limit load: $\lambda^+ \geq \lambda_u$.

\section{\label{sec:conic-repr} Conic representation of strength criteria}
This section deals with the representation of material strength criteria or associated support functions in the framework of conic programming. The first subsection describes the standard conic programming format used by the \texttt{Mosek} solver although other solvers (e.g. CVXOPT, Sedumi, SDPT3) use a format which is quite similar. We then discuss how the conic programming framework is used for limit analysis problems.

\subsection{Conic programming}\label{sec:conic-programming}
Optimization problems entering the \textit{conic programming framework} can be written as:
\begin{equation}
\begin{array}{rl} \displaystyle{\min_\mathbf{x}} & \mathbf{c}\T \mathbf{x} \\
\text{s.t.} &\mathbf{b}_l \leq \mathbf{Ax}\leq \mathbf{b}_u\\
& \mathbf{x}\in \Kk
\end{array} \label{conic-programming}
\end{equation}
where vector $\mathbf{c}$ defines a linear objective functional, matrix $\mathbf{A}$ and vectors $\mathbf{b}_u,\mathbf{b}_l$ define linear inequality (or equality if $\mathbf{b}_u=\mathbf{b}_l$) constraints and where $\Kk=\Kk^1\times \Kk^2\times \ldots \times \Kk^p$ is a product of cones $\Kk^i \subset \RR^{d_i}$ so that $\mathbf{x}\in \Kk \Leftrightarrow \mathbf{x}^i \in \Kk^i \:\forall i=1,\ldots,p$ where $\mathbf{x}=(\mathbf{x}^1,\mathbf{x}^2, \ldots, \mathbf{x}^p)$. These cones can be of different kinds:
\begin{itemize}
\item $\Kk^i = \mathbb{R}^{d_i}$ i.e. no constraint on $\mathbf{x}^i$
\item $\Kk^i = (\mathbb{R}^+)^{d_i}$ is the positive orthant i.e. $\mathbf{x}^i \geq 0$
\item $\Kk^i = \Qq_{d_i}$ the quadratic Lorentz cone defined as:
\begin{equation}
\Qq_{d_i} = \{\mathbf{z} \in \RR^{d_i} \text{ s.t. } \mathbf{z}=(z_0, \bar{\mathbf{z}}) \text{ and } z_0 \geq \|\bar{\mathbf{z}}\|_2\}
\end{equation}
\item $\Kk^i = \Qq_{d_i}^r$ the rotated quadratic Lorentz cone defined as:
\begin{equation}
\Qq_{d_i}^r = \{\mathbf{z} \in \RR^{d_i} \text{ s.t. } \mathbf{z}=(z_0, z_1, \bar{\mathbf{z}}) \text{ and } 2z_0z_1 \geq \|\bar{\mathbf{z}}\|_2^2\}
\end{equation}
\item $\Kk^i = \mathbb{S}_{n_i}^+$ is the cone of semi-definite positive matrices of dimension $n_i$.
\end{itemize}

If $\Kk$ contains only cones of the first two kinds, then the resulting optimization problem \eqref{conic-programming} belongs to the class of Linear Programming (LP) problems. If, in addition, $\Kk$ contains quadratic cones $\Qq_{d_i}$ or $\Qq_{d_i}^r$, then the problem belongs to the class of Second-Order Cone Programming (SOCP) problems. Finally, when cones of the type $\mathbb{S}_{n_i}^+$ are present, the problem belongs to the class of Semi-Definite Programming (SDP) problems. 

\subsection{Conic-representable functions and sets}
In order to use conic programming solvers for our application, limit analysis problems must therefore be formulated following format \eqref{conic-programming}. Inspecting the structure of \eqref{car-stat}, it can be seen that all constraints and the objective functions are linear, except for the strength condition \eqref{eq:crit}. As a result, in order to fit format \eqref{conic-programming}, only \eqref{eq:crit} must be reformulated in a conic sense. Similarly for problem \eqref{car-kin}, only the support function $\pi_G$ must be expressed in conic form to fit the standard format \eqref{conic-programming}. This reformulation step therefore depends on the specific choice of a strength criterion. To do so, we define a generic form of \textit{conic-representable} functions and sets. The \texttt{fenics\_optim} package \cite{jeremy_bleyer_2020_3604086} relies on this specific notion for automating the formulation of generic convex problems. The \texttt{fenics\_optim.limit\_analysis} module accompanying this paper uses these notions and particularizes them for limit analysis problems. 
 
Conic-representable functions are defined as the class of convex functions which can be expressed as follows:
\begin{equation}
\begin{array}{rl} \displaystyle{F(\mathbf{x}) = \min_{\mathbf{y}}} & \mathbf{c}_x\T\mathbf{x}+ \mathbf{c}_y\T\mathbf{y} \\
\text{s.t.} & \mathbf{b}_l\leq \mathbf{Ax}+\mathbf{By} \leq \mathbf{b}_u\\
& \mathbf{y}\in \Kk
\end{array} \label{conic-representable}
\end{equation}
with $\mathbf{x} \in \RR^n$ and in which $\Kk$ is again a product of cones of the kinds detailed in section \ref{sec:conic-programming}. As a by-product of the previous definition, \textit{conic-representable convex sets} correspond to sets for which the indicator function is conic-representable. If $\Kk$ contains only second-order cones (SOC), then we will speak about a \textit{SOC-representable} function. A \textit{SDP-representable} function corresponds to the case when $\Kk$ contains SDP cones, whereas \textit{linear-representable} functions correspond to the case when the conic constraint $\mathbf{y}\in\Kk$ is absent or contains only positive constraints $y_i\geq 0$. 

It is easy to see that if $F$ is SOC-representable (resp. SDP-representable, resp. linear-representable), then its Legendre-Fenchel conjugate $F^*$ is also SOC-representable (resp. SDP-representable, resp. linear-representable). For limit analysis applications, this means that if a strength criterion (expressed as a convex set) is SOC-representable (resp. SDP-representable, resp. linear-representable), then so will be its support function, and vice-versa.

As already mentioned in \cite{bisbos2007second,krabbenhoft2007formulation,makrodimopoulos2010remarks}, a large class of classical strength criteria can be formulated in terms of second-order cone (SOC) constraints or semi-definite positive (SDP) matrix constraints so that their expression or their support function expressions can be expressed in the form \eqref{conic-representable}. Below, we give an example on how conic-representation is used for the case of a 2D plane-strain Mohr-Coulomb material. Many generic conic reformulations can be found in \cite{lobo1998applications, ben2001lectures,boyd2004convex} and especially in the Mosek Modeling Cookbook \cite{mosekcookbook}.

\subsection{Example of the plane-strain Mohr-Coulomb criterion}
In plane-strain conditions, the Mohr-Coulomb criterion with cohesion $c$ and friction angle $\phi$ writes as:
\begin{equation}
\bsig\in \text{MC}_{2D}(c,\phi) \Longleftrightarrow \sqrt{(\sigma_{xx}-\sigma_{yy})^2+4\sigma_{xy}^2} \leq 2c\cos\phi -(\sigma_{xx}+\sigma_{yy})\sin\phi
\end{equation}
which can be also written as:
\begin{align}
&\sqrt{y_1^2+y_2^2} \leq y_0 \Leftrightarrow \mathbf{y}\in \Qq_3\\
\text{with } \begin{Bmatrix}y_0 \\ y_1 \\ y_2 \end{Bmatrix} &= \begin{bmatrix} -\sin\phi & -\sin\phi & 0 \\ 1 & -1 & 0 \\ 0 & 0 & 2 \end{bmatrix}\begin{Bmatrix}\sigma_{xx} \\ \sigma_{yy} \\ \sigma_{xy} \end{Bmatrix} + \begin{Bmatrix}2c\cos\phi \\ 0 \\ 0 \end{Bmatrix}
\end{align}
This expression shows that the criterion is SOC-representable in the sense of format \eqref{conic-representable}.\\

Similarly, its support function is given by:
\begin{align}
\pi(\bd) &= \sup_{\bsig \in \text{MC}_{2D}(c,\phi)} \sigma_{ij}d_{ij} \\
&= \begin{cases}c\cotan\phi \tr(\bd) &\text{if } \tr(\bd) \geq \sin\phi\sqrt{(d_{xx}-d_{yy})^2+4d_{xy}^2} \\ +\infty & \text{otherwise} \end{cases} \notag
\end{align}
which can be expressed as:
\begin{align}
\pi(\bd) = \min_\mathbf{y}\:\: &c\cotan\phi \tr(\bd) \label{pi-MC2D-conic} \\
\text{s.t.} & \begin{bmatrix} 1 & 1 & 0 \\ \sin\phi & -\sin\phi & 0 \\ 0 & 0 & 2\sin\phi \end{bmatrix}\begin{Bmatrix}d_{xx} \\ d_{yy} \\ d_{xy} \end{Bmatrix} = \mathbf{y} \notag \\
&\mathbf{y} \in \Qq_3 \notag
\end{align}
which also fits format \eqref{conic-representable}. Note that the representation is by no means unique since, for instance, we could have equivalently replaced the linear objective term by $c\cotan\phi\, y_0$ or we could also have inverted the linear relation between $\bd$ and $\mathbf{y}$ and introduced $\widetilde{\mathbf{y}}=\mathbf{y}/\sin\phi$ so that:
\begin{align}
\pi(\bd) = \min_{\widetilde{\mathbf{y}}}\:\: &c\cos\phi\, \widetilde{y}_0 \\
\text{s.t.} & \begin{Bmatrix}d_{xx} \\ d_{yy} \\ d_{xy} \end{Bmatrix}=\begin{bmatrix} (\sin\phi)/2 & 1/2 & 0 \\ (\sin\phi)/2 & -1/2 & 0 \\ 0 & 0 & 1/2 \end{bmatrix}\widetilde{\mathbf{y}} \notag \\
&\widetilde{\mathbf{y}} \in \Qq_3 \notag
\end{align}
This expression has the advantage over \eqref{pi-MC2D-conic} of being still well-defined when $\phi=0$, enabling to recover the Tresca/von Mises case.

\subsection{A gallery of conic-representable strength criteria}
The \texttt{fenics\_optim.limit\_analysis} module provides access to a large class of usual strength criteria through the conic formulation of the criterion indicator and support functions, the latter being provided both in terms of strains and velocity discontinuities, see section \ref{sec:continuum-LA}. Table \ref{tab:crit-cont} provides a list of currently available material strength criteria.

\begin{table}
\center
\begin{tabular}{ccc}
\hline 
Strength criterion & Mechanical model & Representation type \\ 
\hline 
\texttt{vonMises} & 2D/3D & SOC\\ 
\hline 
\texttt{DruckerPrager} & 2D/3D & SOC\\ 
\hline 
\texttt{Tresca2D} & 2D & SOC \\ 
\hline 
\texttt{MohrCoulomb2D} & 2D & SOC \\
\hline 
\texttt{Rankine2D} & 2D & SOC \\
\hline
\texttt{Tresca3D} & 3D & SDP \\ 
\hline 
\texttt{MohrCoulomb3D} & 3D & SDP \\ 
\hline 
\texttt{Rankine3D} & 3D & SDP \\ 
\hline 
\texttt{TsaiWu} & 3D & SOC \\ 
\hline 
\end{tabular}
\caption{List of available conic-representable strength criteria for 2D/3D continua} 
\label{tab:crit-cont}
\end{table}
\section{\label{sec:continuum-LA}Finite-element limit analysis of 2D and 3D continua}
In this section, we consider the standard continuum model where $\bSig=\bsig$ is the classical Cauchy symmetric stress tensor, $\bu$ a 2D or 3D velocity field and the associated strain is its symmetric gradient $\bD \bu = \nabla^s \bu$. For simplicity, we will consider the case where there is no fixed loading $\boldsymbol{f}_0=\boldsymbol{t}_0=0$.

\subsection{Kinematic-based formulation}
Let us now consider the finite-element discretization of the kinematic limit analysis approach \eqref{car-kin} for a continuous velocity field and imposed tractions on some part $\partial\Omega_T$ of the boundary:
\begin{equation}
\begin{array}{rl}
\lambda^+ \leq \displaystyle{\inf_{\boldsymbol{u}\in \Vv_h}} & \displaystyle{\int_{\Omega} \pi_G(\nabla^s \boldsymbol u) \dx} \label{kinematic-LA}\\
\text{s.t.} & \int_{\Omega}\boldsymbol{f}\cdot\boldsymbol{u}\dx + \int_{\partial \Omega_T} \boldsymbol{t}\cdot\boldsymbol{u}\dS = 1
\end{array}
\end{equation}
where $\Vv_h$ is a finite-element subspace of Lagrange elements based on a given mesh of typical mesh size $h$. In the above problem, the computed objective function is an upper bound of the exact limit load factor $\lambda^+$ only if the integral of the objective function term is evaluated exactly. In general, this is not possible because of the non-linearity of function $\pi_G$, except in the special case of $\boldsymbol{u}$ being interpolated with $\PP^1$-Lagrange elements so that the gradient is cell-wise constant and the integral becomes trivial.\\

As explained in \cite{makrodimopoulos2007upper}, keeping an exact upper-bound estimate of $\lambda^+$ requires this integral to be estimated by excess. This is possible for a mesh consisting of simplex (straight edges) triangles and a $\PP^2$-Lagrange interpolation for $\boldsymbol{u}$ when using the following so-called \textit{vertex} quadrature scheme:
 \begin{equation}
\int_T F(\mathbf{r}(\boldsymbol{x}))\dx \lesssim \dfrac{|T|}{d+1}\sum_{i=1}^{d+1}  F(\mathbf{r}(\boldsymbol{x}_i)) \label{vertex}
\end{equation}
where $F$ is a convex function and $\mathbf{r}$ is an affinely-varying function over the mesh cell $T$ (either a triangle in dimension $d=2$ or a tetrahedron for $d=3$ of area/volume $|T|$) with $\mathbf{r}(\boldsymbol{x}_i)$ being its value at the $d+1$-vertices. In this case, we have:
\begin{equation}
\begin{array}{rl}
\lambda^+ \leq \lambda_u = \displaystyle{\inf_{\boldsymbol{u}\in \Vv_h}} & \displaystyle{\sum_{T\in \Tt_h}\dfrac{|T|}{d+1} \sum_{i=1}^{d+1} \pi_G(\nabla^s \boldsymbol u(\boldsymbol{x}_i))} \\
\text{s.t.} & \displaystyle{\int_{\Omega}\boldsymbol{f}\cdot\boldsymbol{u}\dx + \int_{\partial \Omega_T} \boldsymbol{t}\cdot\boldsymbol{u}\dS = 1}
\end{array}
\end{equation}

As discussed in \cite{bleyer2019automating}, the \texttt{fenics\_optim} package enables to solve convex variational problem of the form:
\begin{equation}
\begin{array}{rl} \displaystyle{\inf_{\boldsymbol{u}\in \Vv}} & \displaystyle{\int_{\Omega} \left(j_1\circ \ell_1 (\boldsymbol{u})+\ldots+j_p\circ \ell_p(\boldsymbol{u})\right)\dx} \\
\text{s.t.}  &\boldsymbol{u}\in \Kk
\end{array} \label{co-prob}
\end{equation}
where $j_i$ are conic-representable convex functions and $\ell_i$ are linear operators which can be expressed using UFL symbolic operators. Each individual term $j_i$ can be specified independently along with its conic representation and the quadrature scheme used for the computation over $\Omega$ so that \eqref{co-prob} is in fact of the form:
\begin{equation}
\begin{array}{rl} \displaystyle{\inf_{\boldsymbol{u}\in \Vv}} & \displaystyle{\sum_{i=1}^p\sum_{g_i = 1}^{N_{g,i}} \omega_{g_i}j_i\left(\ell_i (\boldsymbol{u};\boldsymbol{x}_{g_i})\right)} \\
\text{s.t.}  &\boldsymbol{u}\in \Kk
\end{array} \label{co-prob-discr}
\end{equation} 
where $\ell_i (\boldsymbol{u};\boldsymbol{x}_{g_i})$ denotes the evaluation of $\ell_i (\boldsymbol{u})$ at a quadrature point $\bx_{g_i}$.

Local auxiliary variables $\mathbf{y}$ of the conic representation \eqref{conic-representable} for each $j_i$ will be added to the optimization problem for each quadrature point $\boldsymbol{x}_{g_i}$. More details can be found in \cite{bleyer2019automating}.\\

Obviously, problem \eqref{kinematic-LA} is a problem of the form \eqref{co-prob} in which $j_1 = \pi_G$, $\ell_1 = \nabla^s$ and $j_2$ is the indicator function of the linear constraint \mbox{$\int_{\Omega}\boldsymbol{f}\cdot\boldsymbol{u}\dx + \int_{\partial \Omega_T} \boldsymbol{t}\cdot\boldsymbol{u}\dS = 1$} with $\ell_2$ being the identity. Ignoring import statements and mesh generation, we now give a few lines of Python script to illustrate how the \texttt{fenics\_optim.limit\_analysis} module enables to formulate an upper bound limit analysis problem for a 2D Mohr-Coulomb material using a $\PP^2$-Lagrange interpolation for $\boldsymbol{u} \in \Vv_h$. First, the corresponding function space \texttt{V} is defined and fixed boundary conditions are imposed on the part named "\texttt{border}" of the boundary. A \texttt{MosekProblem} object is instantiated and a first optimization field \texttt{u} belonging to function space \texttt{V} is added to the problem and is constrained to satisfy the Dirichlet boundary conditions:
\begin{pythoncode}
V = VectorFunctionSpace(mesh, "CG", 2)
bc = DirichletBC(V, Constant((0.,0.)), border)

prob = MosekProblem("Upper bound limit analysis")
u = prob.add_var(V, bc=bc)
\end{pythoncode}
The external work normalization constraint is then added by defining the function space for the Lagrange multiplier corresponding to the constraint (here it is scalar so we use a \texttt{"Real"} function space) and passing the corresponding constraint in its weak form as follows (here $\boldsymbol{t}=0$):
\begin{pythoncode}
R = FunctionSpace(mesh, "R", 0)
def Pext(lamb):
    return [lamb*dot(f,u)*dx]
prob.add_eq_constraint(R, A=Pext, b=1)
\end{pythoncode}
Now, a Mohr-Coulomb material is instantiated and provides access to its convex support function. The input arguments are the strain $\nabla^s \boldsymbol{u}$ written in terms of UFL operators as well as the choice for the quadrature scheme used for its numerical evaluation. Here the \textit{vertex} scheme \eqref{vertex} is chosen. This convex function is then added to the problem before asking for its optimization by \texttt{Mosek}.
\begin{pythoncode}
mat = MohrCoulomb2D(c, phi)
strain = sym(grad(u))
pi = mat.support_function(strain, quadrature_scheme="vertex")
prob.add_convex_term(pi)

prob.optimize()
\end{pythoncode}

\subsection{\label{sec:static-continuum}Static-based formulation}
Conversely, a lower bound estimate of the limit load factor $\lambda^+$ can be obtained by considering a statically admissible discretization of the stress field $\bsig$. For this purpose, we consider the lower bound element of \cite{lysmer1970limit,sloan1988lower}, which includes statically admissible discontinuities between facets. The problem we aim at solving is the following:
\begin{equation}
\begin{array}{rl}
\displaystyle{\lambda^+ \geq \lambda_s = \sup_{\lambda\in \RR, \bsig \in \Ww_h}} & \lambda \\
\text{s.t. }\: & \div\bsig + \lambda \boldsymbol{f} = 0 \quad\text{in } \Omega\\
& \jump{\bsig}\cdot\boldsymbol{n} = 0 \quad\text{through } \Gamma\\
& \bsig\cdot \boldsymbol{n} = \lambda \boldsymbol{t} \quad\text{on } \partial\Omega_T\\
& \bsig(\boldsymbol{x}) \in G \quad \forall \boldsymbol{x} \in \Omega
\end{array} \label{static-LA}
\end{equation}
where $\Gamma$ is the set of internal facets and $\jump{\bsig}$ the stress discontinuity across this facet of normal $\boldsymbol{n}$.\\

For this problem, we have as optimization variables one real $\lambda$ and a discontinuous piecewise-linear field $\bsig$ represented by a vector of dimension 3 in 2D (6 in 3D), $\Ww_h$ being the associated discontinuous $\PP^1_d$ function space:
\begin{pythoncode}
prob = MosekProblem("Lower bound limit analysis")
R = FunctionSpace(mesh, "R", 0)
W = VectorFunctionSpace(mesh, "DG", 1, dim=3)
lamb, Sig = prob.add_var([R, W])
sig = as_matrix([[Sig[0], Sig[2]],
                 [Sig[2], Sig[1]]])
\end{pythoncode}

Assuming piecewise constant body force over each cell, the first local equilibrium equation can be equivalently written weakly using $\PP^0$ velocity fields as Lagrange multipliers:
\begin{pythoncode}
V_eq = VectorFunctionSpace(mesh, "DG", 0)
def equilibrium(u):
    return [dot(u,f)*lamb*dx, dot(u,div(sig))*dx]
prob.add_eq_constraint(V_eq, A=equilibrium)
\end{pythoncode}
Note that in the \texttt{equilibrium} function each block respectively corresponds to the optimization variable $\lambda$ then $\bsig$, in the order they have been initially defined. See more details in \cite{bleyer2019automating} on the underlying block-structure of the problem.

Similarly, discontinuous affine displacements $\boldsymbol{v}$ defined on the mesh facets only\footnote{They are called \texttt{Discontinuous Lagrange Trace} elements in FEniCS} are used as Lagrange multipliers for the second and third constraints (here $\boldsymbol{t}=0$):
\begin{pythoncode}
V_jump = VectorFunctionSpace(mesh, "Discontinuous Lagrange Trace", 1)
def continuity(v):
    return [None, dot(avg(v),jump(sig,n))*dS
                + dot(v, dot(sig,n))*ds(0)]
prob.add_eq_constraint(V_jump, A=continuity)
\end{pythoncode}
where the \texttt{dS} term corresponds to the integral over all internal facets $\Gamma$ and \texttt{ds(0)} corresponds to the integral over $\partial \Omega_T$ in this case.

Finally, the problem objective function is added to the problem as well as the strength criterion condition. The latter is treated as a convex function through its indicator function. A quadrature scheme is still needed since quadrature points will correspond to points at which the strength condition will be enforced. In the present case, the \textit{vertex} scheme will enforce the strength condition at the triangle vertices so that it will be satisfied everywhere inside the cell by convexity. The maximization problem is then solved by \texttt{Mosek}:
\begin{pythoncode}
prob.add_obj_func([1, None])

crit = mat.criterion(Sig, quadrature_scheme="vertex")
prob.add_convex_term(crit)

prob.optimize(sense="maximize")
\end{pythoncode}

Note that \texttt{Mosek} solutions give access to dual variables (Lagrange multipliers) so that a pseudo velocity field can be obtained from \texttt{u} for instance.

\subsection{\label{sec:mixed-FE}Mixed finite-element discretizations}
The kinematic formulation \eqref{kinematic-LA} is sometimes difficult to use because the support function expression may be cumbersome to derive. As a result, a static-based formulation as in \eqref{static-LA} is usually more attractive as it requires only the expression of the strength criterion itself. Formulation \eqref{car-mixed} enables to produce upper bounds equivalent to \eqref{kinematic-LA}, provided a proper choice of interpolation spaces for $\bu$ and $\bSig$ and quadrature rules. In the general case, such a formulation forms the basis of mixed finite-element interpolations for which the bounding status is lost in general \cite{capurso1971limit,anderheggen1972finite,casciaro1982mixed,christiansen1999computation}. This formulation reads here as:
\begin{equation}
\begin{array}{rl}
\displaystyle{\lambda_m = \sup_{\lambda\in \RR, \bsig \in \Ww_h}} & \lambda \\
\text{s.t. }\: & \displaystyle{\int_{\Omega}\bsig:\nabla^s\boldsymbol{u}\dx = \lambda \left(\int_{\Omega}\boldsymbol{f}\cdot\boldsymbol{u}\dx + \int_{\partial\Omega_T} \boldsymbol{t}\cdot\bu \dS\right)} \quad\forall \boldsymbol{u}\in \Vv_h\\
& \bsig(\boldsymbol{x}) \in G \quad \forall \boldsymbol{x} \in \Omega
\end{array} \label{mixed-LA}
\end{equation}
in which the static equilibrium conditions have been replaced by their weak counterpart using the virtual work principle for a class of kinematically admissible continuous velocity fields $\boldsymbol{u}\in \Vv_h$. As mentioned in \cite{bleyer2019automating}, if $\Vv_h$ corresponds to a continuous $\PP^1$-Lagrange interpolation, $\Ww_h$ to a discontinuous $\PP^0$ interpolation of the stress field and the strength criterion is enforced at one point in each cell (this is enough since $\bsig$ is cell-wise constant), then \eqref{mixed-LA} is equivalent, it is even the dual problem, to \eqref{kinematic-LA} for the same $\Vv_h$. 

For other cases, such as $\Vv_h$ being $\PP^2$-Lagrange and $\Ww_h$ discontinuous $\PP^1_d$-Lagrange, quadrature rules must be specified both for the equilibrium constraint as well as for the strength criterion (criterion enforcement points). If both quadrature rules are identical, then problem \eqref{mixed-LA} is equivalent to problem $\eqref{kinematic-LA}$ for the chosen quadrature rule. Upon choosing a vertex scheme, the objective value $\lambda_m$ will be an upper-bound to $\lambda^+$. Other choices such as quadrature points located at the mid-side points do not produce rigorous upper bounds but are usually observed to converge from above in practice. When both quadrature rules are different, one obtains a true mixed-interpolation and again the bounding status is lost.

The FEniCS formulation for the first case $\PP^1/\PP^0$ for $\Vv_h/\Ww_h$ would read as (here quadrature rules are trivial one-point rules by default):
\begin{pythoncode}
prob = MosekProblem("Upper bound from static formulation")
R = FunctionSpace(mesh, "R", 0)
W = VectorFunctionSpace(mesh, "DG", 0, dim=3)
lamb, Sig = prob.add_var([R, W])

V = VectorFunctionSpace(mesh, "CG", 1)
bc = DirichletBC(V, Constant((0, 0)), border)

sig = as_matrix([[Sig[0], Sig[2]],
                 [Sig[2], Sig[1]]])
def equilibrium(u):
    return [lamb*dot(u, f)*dx, -inner(sig, sym(grad(u)))*dx]
prob.add_eq_constraint(V, A=equilibrium, bc=bc)

prob.add_obj_func([1, None])

crit = mat.criterion(Sig)
prob.add_convex_term(crit)

prob.optimize(sense="maximize")
\end{pythoncode}

\subsection{Discontinuous velocity fields}
As mentioned in the introduction, the use of discontinuous velocity fields is interesting in a limit analysis context due to a higher accuracy and the absence of volumetric or shear locking effects.

In case of discontinuous velocity fields across a set $\Gamma$ of discontinuity surfaces, the kinematic limit analysis formulation \eqref{kinematic-LA} now becomes:
\begin{equation}
\begin{array}{rl}
\lambda^+ \leq \lambda_u = \displaystyle{\inf_{\boldsymbol{u}\in \Vv_h}} & \displaystyle{\int_{\Omega} \pi_G(\nabla^s \boldsymbol u) \dx + \int_{\Gamma} \Pi_G(\jump{\boldsymbol u};\boldsymbol{n}) \dS } \label{kinematic-LA-disc}\\
\text{s.t.} & \int_{\Omega}\boldsymbol{f}\cdot\boldsymbol{u}\dx + \int_{\partial \Omega_T} \boldsymbol{t}\cdot\boldsymbol{u}\dS = 1
\end{array}
\end{equation}
where the second term of the objective function denotes the dissipated power contribution of the velocity jumps $\jump{\bu}$ through a surface $\Gamma$ of normal $\bn$. It is computed from the discontinuity support function: 
\begin{equation}
\Pi_G(\bv;\bn) = \sup_{\bsig \in G} \{(\bsig\cdot\bn)\cdot \bv\} = \pi_G\left(\bv \overset{s}{\otimes} \bn\right)
\end{equation}
i.e. the jump operator is here $\bJ\bu = \jump{\bu}\overset{s}{\otimes} \bn$. As a result, the second term will also be conic-representable and will be treated similarly. Despite this relation, the conic representation of $\Pi_G(\bv;\bn)$ is usually implemented explicitly in the material library due to the potential savings in terms of auxiliary variables compared to calling directly $\pi_G\left(\bv \overset{s}{\otimes} \bn\right)$. A local projection of $\jump{\bu}$ on the facet $(\bn, \boldsymbol{t})$-plane is performed as $\Pi_G$ functions are in general naturally expressed in this local frame. 
Finally, adding the discontinuity term to the optimization problem is similar to the first term, one must just specify the integration measure (the inner facets) and the quadrature rule which must be used to perform the facet integration.

\subsection{Mesh refinement}
For improving the quality of the computed limit load estimates, a mesh refinement procedure is also implemented based on the contribution of each cell to maximum resisting work $\Pp^{(mr)}$ in the context of a kinematic approach. The cell contributions are sorted in descending order and we compute the cumulated contribution to $\Pp^{(mr)}$ until reaching a user-specified threshold $\eta \Pp^{(mr)}$ (with a ratio $\eta$ of 0.5 typically) of the total dissipation. The first $k$ cells whose total contribution is at least $\eta \Pp^{(mr)}$ are then marked for mesh refinement. In the case of discontinuous elements, the facet contribution to $\Pp^{(mr)}$ is computed for each facet, split evenly between the two sharing cells and added to the cell contributions. For the lower bound approach, a similar procedure is used except that we use the dual variable associated with the local equilibrium equation to reconstruct a piecewise linear velocity field from which we compute a cell-wise contribution to the total dissipation.

\subsection{Vertical cut-off stability}
\begin{figure}
\begin{center}
\includegraphics[width=0.57\textwidth]{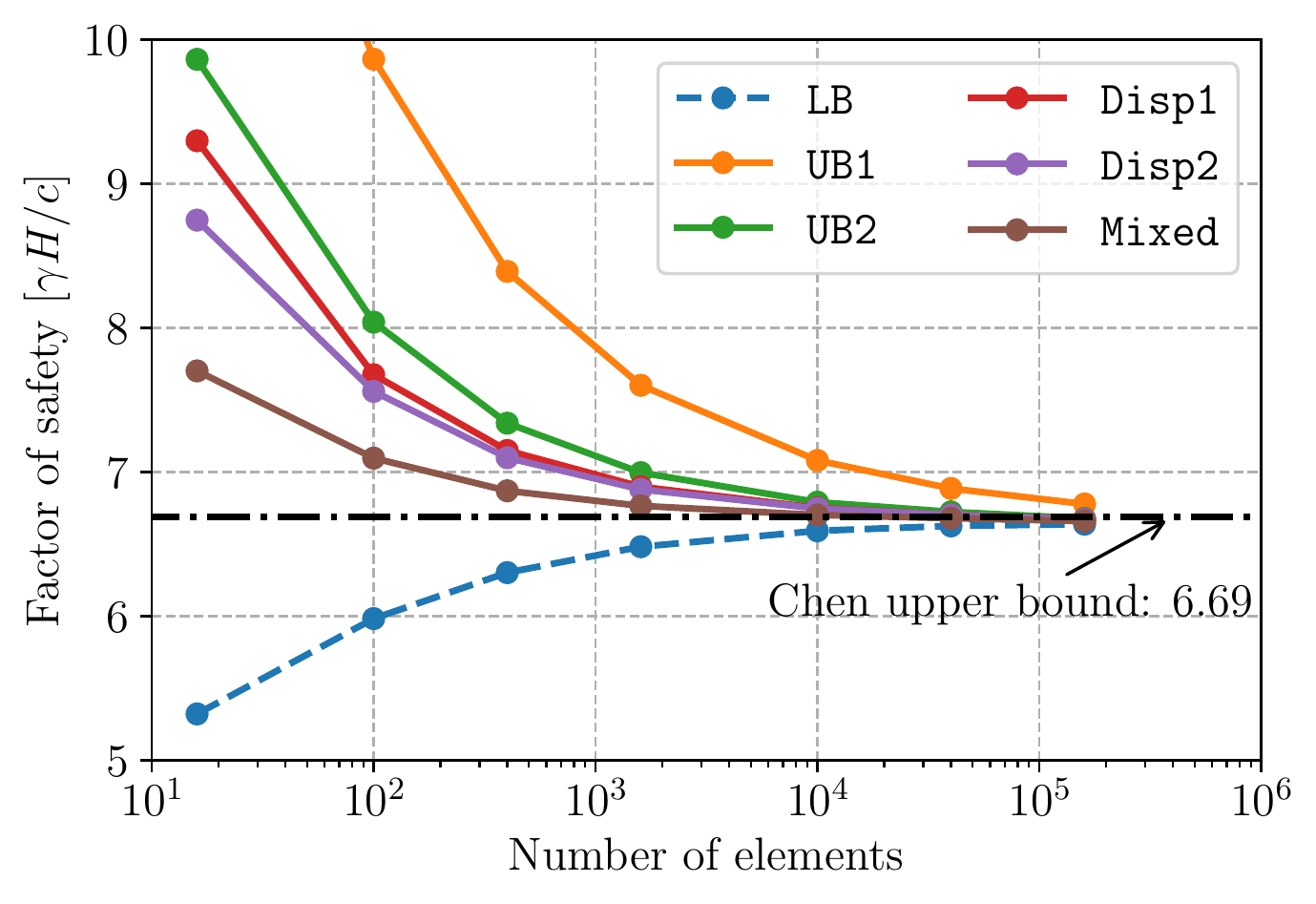}
\includegraphics[width=0.42\textwidth]{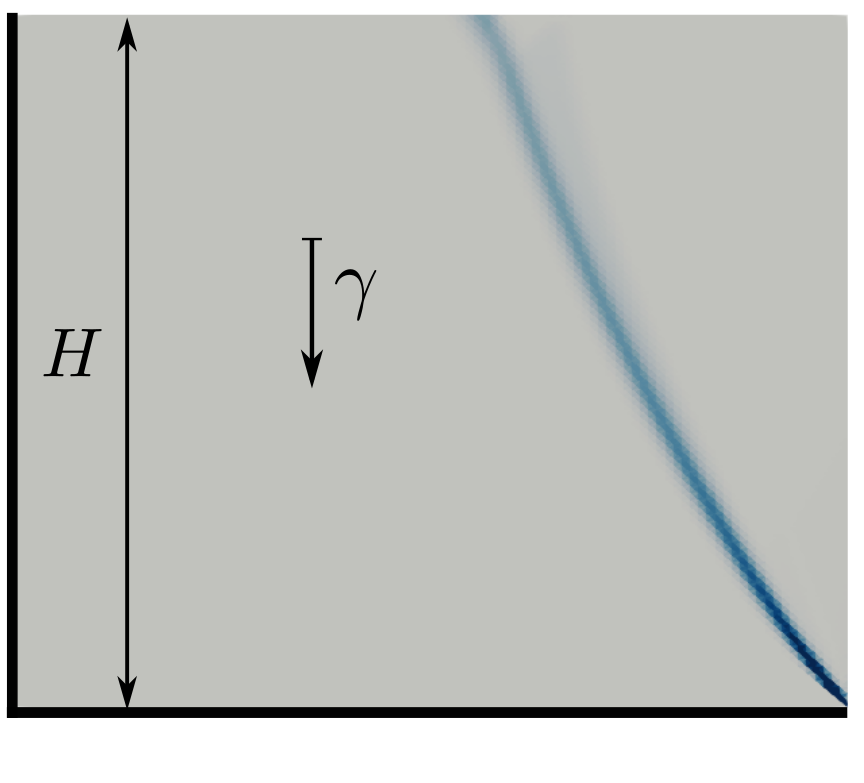}
\end{center}
\caption{Left: Convergence of the vertical slope factor of safety for various finite element discretization: \texttt{UB1} (resp. \texttt{UB2}) are the $\PP^1$ (resp. $\PP^2$) Lagrange upper bound elements, \texttt{LB} the lower bound element and \texttt{Disp1}, \texttt{Disp2}, \texttt{Mixed} correspond to non-upper bound elements considered in \cite{krabbenhoft2007formulation}. Analytical upper bound from Chen \cite{chen2013limit}. Right: local dissipation map}
\label{slope-results}
\end{figure}

\begin{figure}
\begin{center}
\includegraphics[width=0.6\textwidth]{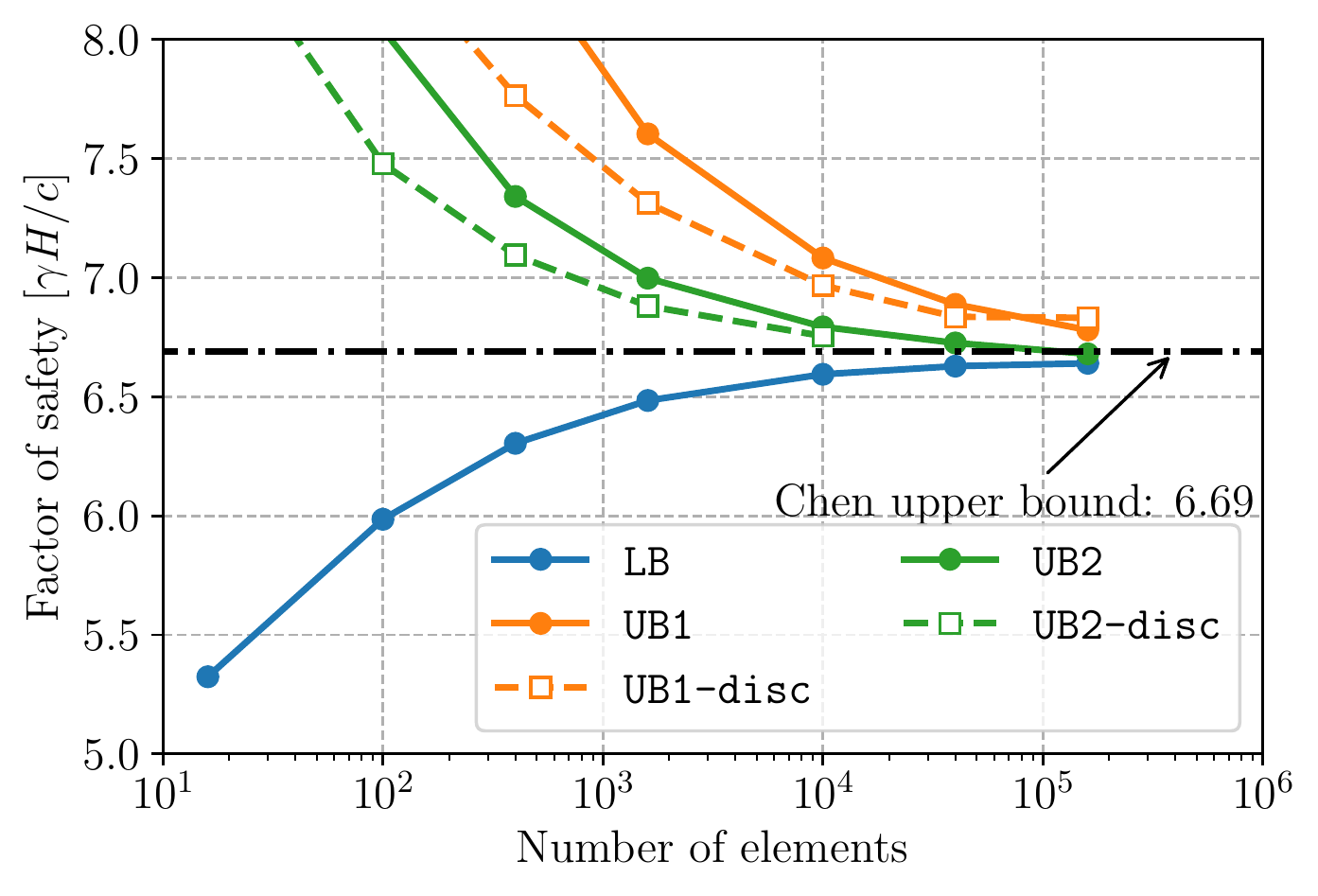}
\end{center}
\caption{Convergence of the vertical slope factor of safety for various finite element discretization for continuous and discontinuous upper bound elements. Analytical upper bound from Chen \cite{chen2013limit}.}
\label{slope-results-disc}
\end{figure}

\begin{figure}
\begin{center}
\includegraphics[width=0.55\textwidth]{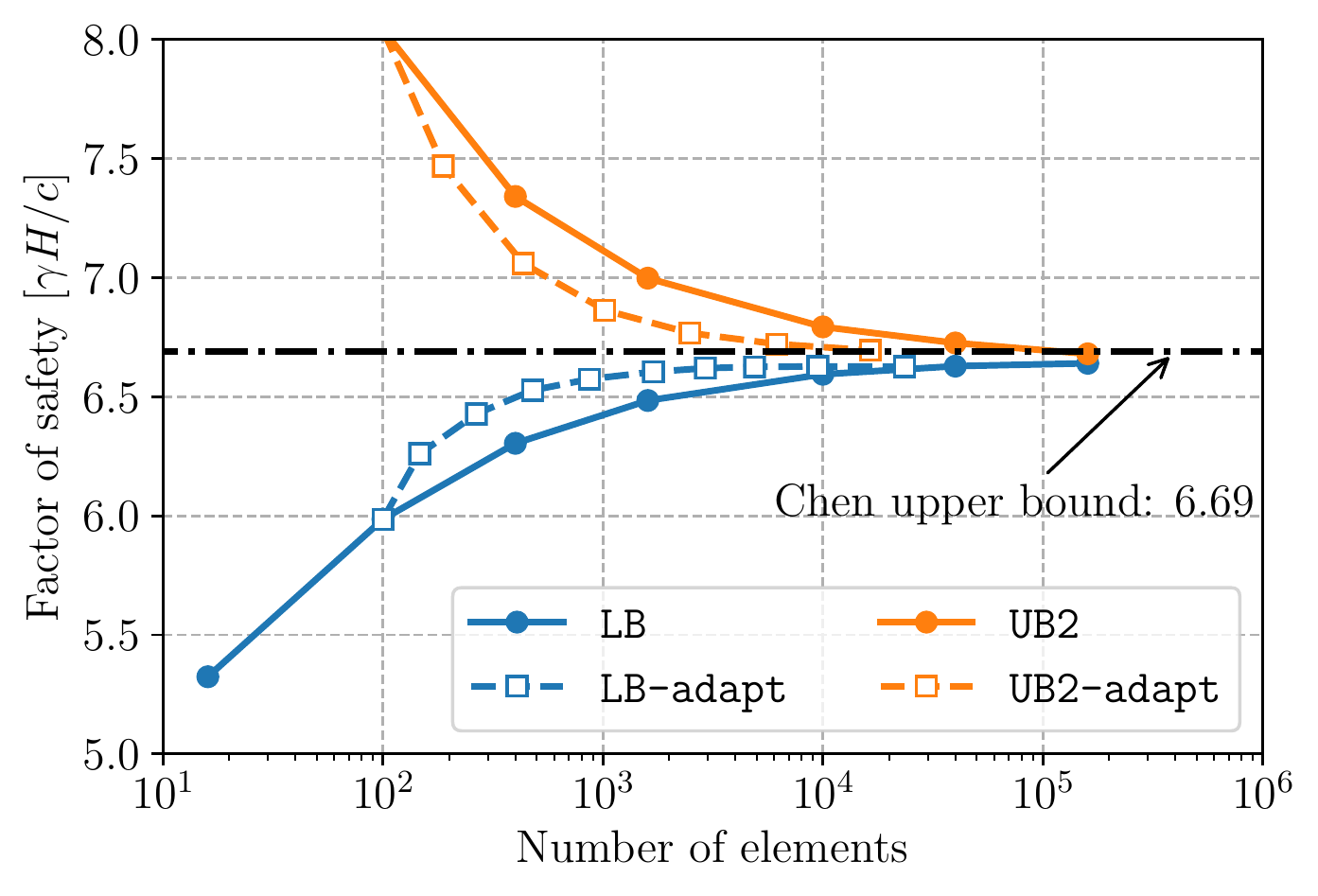}
\includegraphics[width=0.44\textwidth]{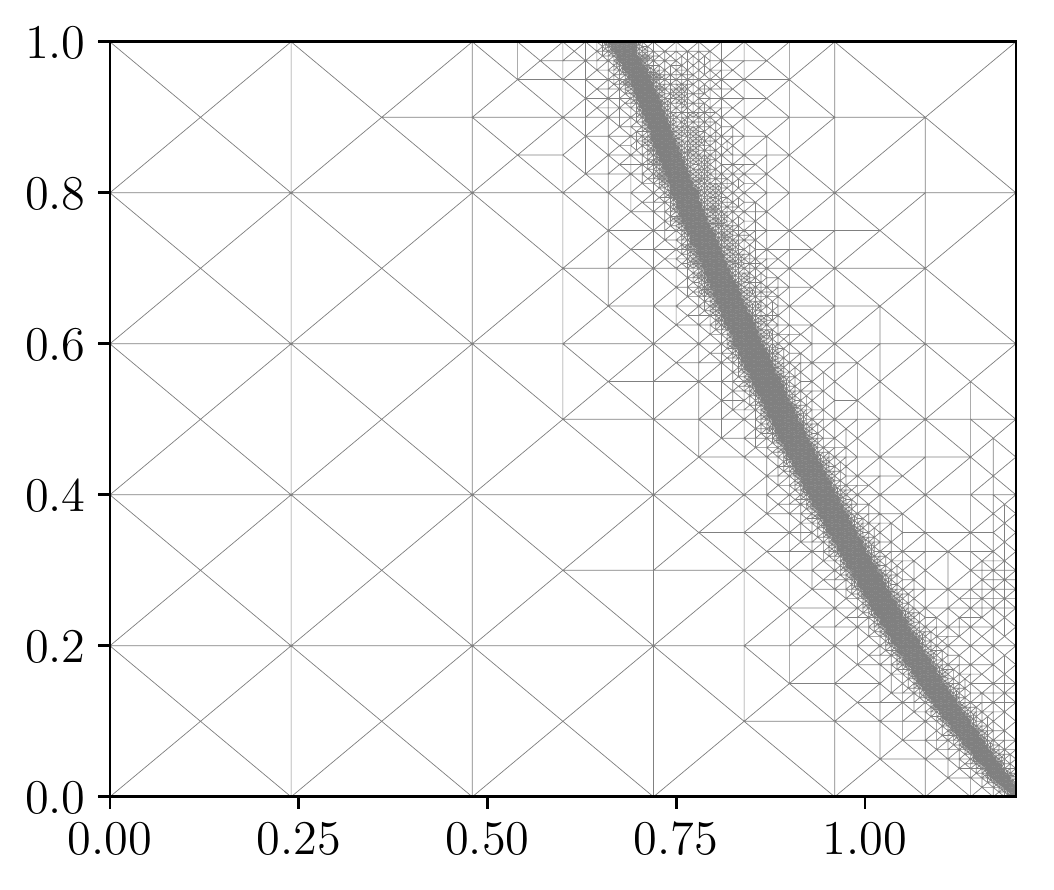}
\end{center}
\caption{Left: Convergence of the vertical slope factor of safety for adapted meshes. Analytical upper bound from Chen \cite{chen2013limit}. Right: adapted mesh after 7 refinement steps}
\label{slope-results-adapt}
\end{figure}

The different discretization choices are illustrated on the stability of a vertical slope under its self-weight $\gamma$. The rectangular domain, of height $H$, is clamped on both bottom and left boundaries and traction-free on the remaining boundaries (see Figure \ref{slope-results}-right). The soil is modelled as a Mohr-Coulomb material with cohesion $c$ and friction angle $\phi=30^{\circ}$ under plane strain conditions. The slope factor of safety is given by the maximum value of the non-dimensional quantity $(\gamma H/c)^+$. Convergence of he factor of safety estimates for various finite-element discretizations are reported on Figure \ref{slope-results}-left as a function of the total number of elements, the concentration of the local dissipation $\pi_G(\nabla^s\bu)$ along the slip-line is represented on Figure \ref{slope-results}-right. The comparison between continuous and discontinuous upper bound finite elements is reported in Figure \ref{slope-results-disc}. Similarly, the comparison between uniform and adaptive mesh refinement is reported in Figure \ref{slope-results-adapt} along with the final adapted mesh.

\section{\label{sec:homog}Computation of macroscopic strength properties through homogenization theory}
In this section, we show how to adapt the above-described framework to compute macroscopic strength properties of heterogeneous materials. Indeed, limit analysis can be combined to homogenization theory \cite{suquet1985elements,de1986fundamental} to compute macroscopic strength criteria of heterogeneous materials by solving an auxiliary limit analysis problem formulated on a representative elementary volume or a periodic unit-cell for periodically heterogeneous materials. We will consider here the latter case and denote $\Aa$ the unit cell domain. The strength conditions are described by a spatially varying strength criterion $G(\boldsymbol{x})$. The main goal of homogenization in limit analysis is to compute the macroscopic strength domain $G^{hom}$ defined as follows:
\begin{equation}
\bSig \in G^{hom} \Longleftrightarrow \begin{cases}
\div \bsig = 0 &\text{in } \Aa\\
\jump{\bsig}\cdot\bn = 0 &\text{through } \Gamma\\
\bsig\cdot\bn & \text{antiperiodic} \\
\bsig(\boldsymbol{x}) \in G(\boldsymbol{x}) & \forall \boldsymbol{x}\in \Aa\\
\displaystyle{\dfrac{1}{|\Aa|}\int_\Aa \bsig \dx = \bSig}
 \end{cases} \label{Ghom}
\end{equation}

Upon choosing a given loading direction $\bSig_0$ of arbitrary magnitude, one can look for an estimate to the maximum load factor $\lambda^+$ such that $\lambda^+ \bSig_0 \in G^{hom}$, yielding the following maximization problem:

\begin{equation}
\begin{array}{rll}
\displaystyle{\lambda^+ \geq \lambda_s= \sup_{\lambda\in \RR, \bsig \in \Ww_h}} & \lambda & \\
\text{s.t. }\: & \div\bsig = 0 &\text{in } \Aa\\
& \jump{\bsig}\cdot\boldsymbol{n} = 0 &\text{through } \Gamma\\
& \bsig\cdot\bn & \text{antiperiodic} \\
& \bsig(\boldsymbol{x}) \in G(\boldsymbol{x}) & \forall \boldsymbol{x} \in \Aa\\
& \displaystyle{\dfrac{1}{|\Aa|}\int_\Aa \bsig \dx = \lambda\bSig_0}
\end{array} \label{static-Ghom}
\end{equation}

Definition \eqref{static-Ghom} correspond to the static formulation of a limit analysis problem formulated on $\Aa$ with loading being parametrized by $\bSig_0$. Its dual counter-part (kinematic formulation) can be shown to be given by\footnote{This formulation can obviously be easily extended to take into account discontinuous velocity fields.}:
\begin{equation}
\begin{array}{rl}
\lambda^+ \leq \lambda_u = \displaystyle{\inf_{\boldsymbol{D}\in\RR^6, \boldsymbol{u}\in \Vv_h}} & \displaystyle{\int_{\Omega} \pi(\boldsymbol{D}+\nabla^s \boldsymbol u;\boldsymbol{x}) \dx} \label{kinematic-Ghom}\\
\text{s.t.} & |\Aa|\bSig_0:\boldsymbol{D} = 1 \\
& \boldsymbol{u} \text{ periodic}
\end{array}
\end{equation}
 Let us remark that:
\begin{equation}
\Pi_{hom}(\boldsymbol{D}) = \lambda^+\bSig_0:\boldsymbol{D} \leq \dfrac{1}{|\Aa|}\int_{\Omega} \pi(\boldsymbol{D}+\nabla^s \boldsymbol u;\boldsymbol{x}) \dx 
\end{equation}
where $\displaystyle{\Pi_{hom}(\boldsymbol{D}) := \sup_{\bSig \in G^{hom}} \{\bSig:\boldsymbol{D}\}}$ is the macroscopic support function.

In \eqref{kinematic-Ghom}, the macroscopic strain is considered as an additional optimization variable since the loading direction $\bSig_0$ is fixed. It is also possible to prescribe directly the macroscopic strain direction $\boldsymbol{D}$, leaving free the loading direction, by removing the first constraint in \eqref{kinematic-Ghom}.\\

\begin{figure}
\begin{center}
\begin{subfigure}[b]{0.49\textwidth}
\includegraphics[width=0.9\textwidth]{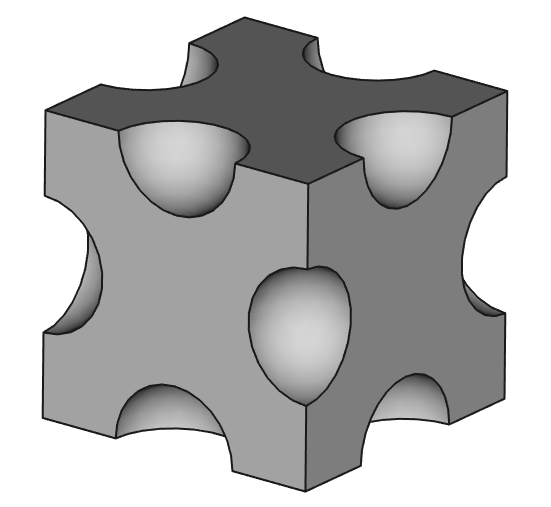}
\caption{Porous 3D medium unit cell}
\label{porous-UC}
\end{subfigure}
\hfill
\begin{subfigure}[b]{0.49\textwidth}
\includegraphics[width=\textwidth]{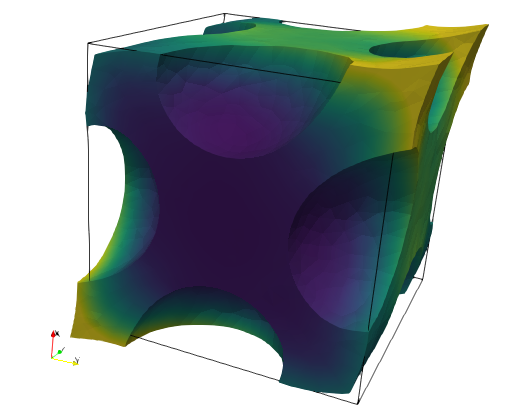}
\caption{Unit-cell response to a pure-shear loading}
\label{porous-shear}
\end{subfigure}
\end{center}
\caption{Homogenization of a 3D porous medium}
\end{figure}

\begin{figure}
\begin{center}
\begin{subfigure}[b]{0.39\textwidth}
\includegraphics[width=\textwidth]{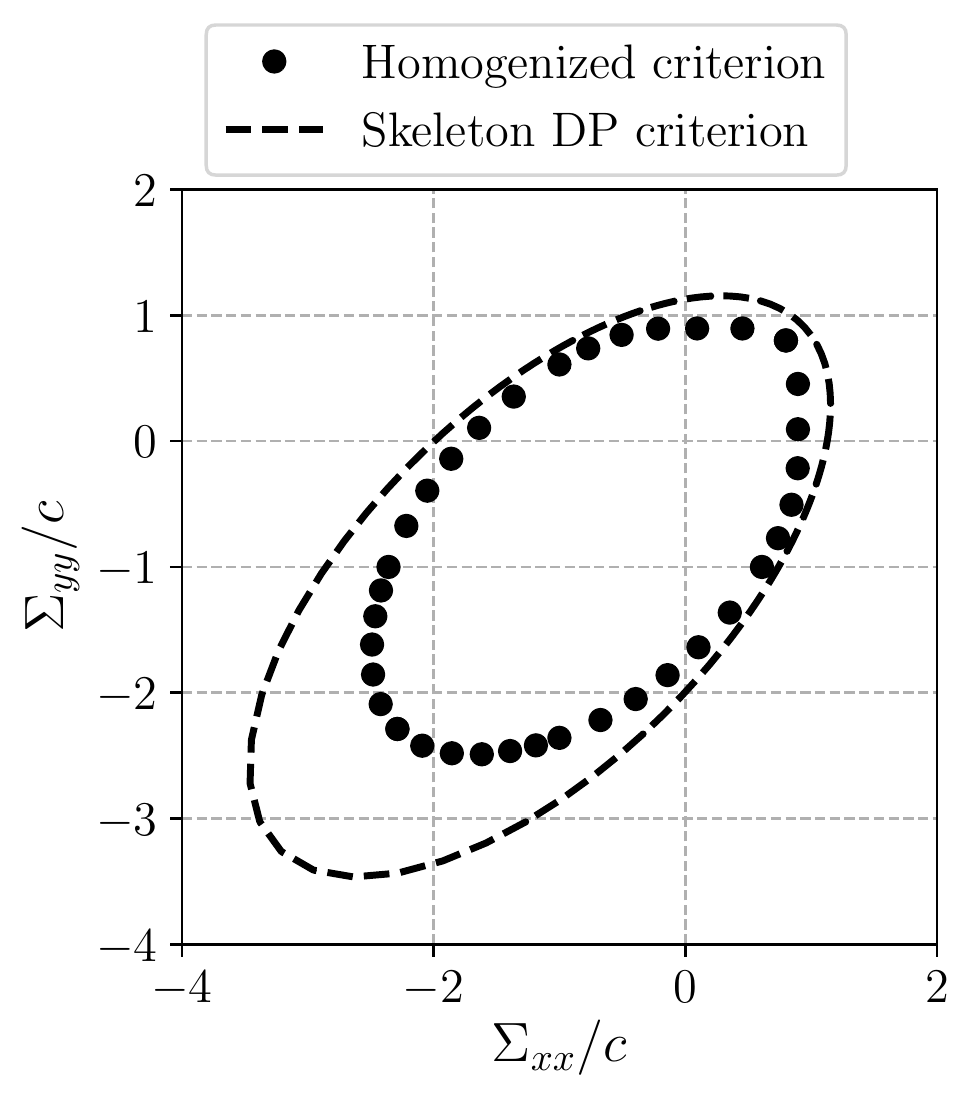}
\caption{$(\bSig_{xx},\bSig_{yy})$-plane}
\end{subfigure}
\hfill
\begin{subfigure}[b]{0.6\textwidth}
\includegraphics[width=\textwidth]{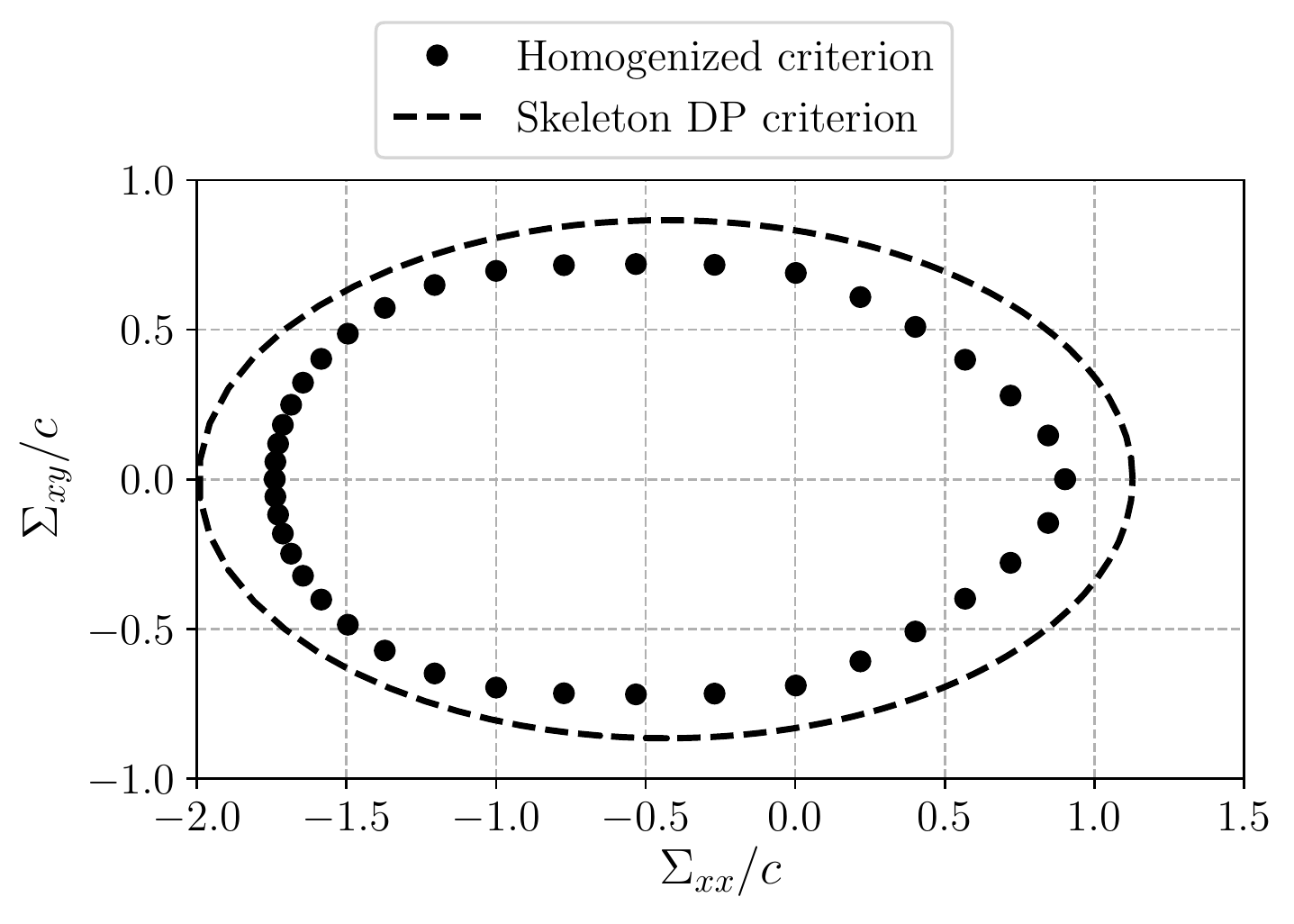}
\caption{$(\bSig_{xx},\bSig_{xy})$-plane}
\end{subfigure}
\end{center}
\caption{Homogenized strength criterion}
\label{porous-Ghom}
\end{figure}

Formulation \eqref{kinematic-Ghom} is applied to a 3D periodic porous medium with pores following a face-centered cubic system (see Figure \ref{porous-UC}) made of a Drucker-Prager material $(c=1,\phi=30^{\circ})$. The unit cell response $\boldsymbol{U}(\boldsymbol{x}) = \boldsymbol{D}\cdot \boldsymbol{x}+\boldsymbol{u}(\boldsymbol{x})$ to a pure shear loading $\bSig_{xy,0}=1$ is represented on Figure \ref{porous-shear}. The macroscopic strength domain in the $(\bSig_{xx},\bSig_{yy})$ plane and in the $(\bSig_{xx},\bSig_{xy})$ with other $\bSig_{ij}=0$ have been represented on Figure \ref{porous-Ghom}. It can be observed how much the presence of the pores reduces the original Drucker-Prager criterion of the skeleton.

\section{\label{sec:plates}Plates and shells}
\subsection{Thin and thick plates}
In this section, we show how limit analysis of thin and thick plates can be tackled by following exactly the general format described in section \ref{sec:generalized-LA} and implemented easily by taking advantage of the high-level symbolic representation of operators in FEniCS. First, let us recall that the kinematic limit analysis of thin plates obeying a Love-Kirchhoff kinematics corresponds to the following minimization problem:
\begin{equation}
\begin{array}{rl}
\displaystyle{\inf_{w}} & \displaystyle{\int_\Omega \pi_G(\nabla^2 w)\dx +\int_\Gamma \Pi_G(\jump{\partial_n w} \boldsymbol{n} ; \boldsymbol{n})\dS} \\
\text{s.t.} & \displaystyle{\int_{\Omega} fw \dx = 1}
\end{array} \label{BENDING-PLATE}
\end{equation}
for a distributed loading $f$. It can be seen that \eqref{BENDING-PLATE} complies with \eqref{car-kin} where the generalized velocity is only the virtual deflection $w$, the strain operator is the curvature $\bD\bu = \nabla^2 w$ and the velocity jump is given by a normal rotation jump $\bJ\bu = \jump{\partial_n w}\boldsymbol{n}\otimes\boldsymbol{n}$. In the above, $\pi_G$ and $\Pi_G$ are the support functions of the corresponding thin plate strength criterion $G_\text{bend}$, expressed solely on the bending moment tensor $\boldsymbol{M}=\begin{bmatrix}
M_{xx} & M_{xy}\\ M_{xy} & M_{yy}
\end{bmatrix}$ with $\Pi_G(\jump{\partial_n w} \boldsymbol{n} ; \boldsymbol{n}) = \pi_G( \jump{\partial_n w}\boldsymbol{n}\otimes\boldsymbol{n})$. As a result, thin plate strength criteria and their support functions are treated exactly as 2D continua with the only difference coming from the definition of the strain and discontinuity operators. We refer to \cite{bleyer2019automating} for an implementation example of thin plate limit analysis using the \texttt{fenics\_optim} package. Let us also remark that implementing the corresponding lower bound static approach is more involved due to the more complicated continuity conditions involving equivalent Kirchhoff shear forces (see \cite{krabbenhoft2002lower,bleyer2014computational}).

As regards thick plates involving shear and bending strength conditions and following Reissner-Mindlin kinematics, the main difference with respect to 2D/3D continua or thin plates is that two unknown fields must be considered instead of one: bending moments $\boldsymbol{M}$ and shear forces $\boldsymbol{Q}=(Q_x,Q_y)$ for the static approach and out-of-plane deflection $w$ and plate rotations $\boldsymbol{\theta}=(\theta_x,\theta_y)$ for the kinematic approach. We now discuss only the latter since the former will be discussed in the more general case of shells in the next subsection. Following \cite{bleyer2015kin_thick_plates}, kinematic limit analysis of thick plates can be written as:
\begin{equation}
\begin{array}{rl}
\displaystyle{\inf_{w, \btheta}} & \displaystyle{\int_\Omega \pi_G((\nabla \btheta, \nabla w-\btheta))\dx +\int_\Gamma \Pi_G((\jump{\btheta},\jump{w})  ; \boldsymbol{n})\dS} \\
\text{s.t.} & \displaystyle{\int_{\Omega} fw \dx = 1}
\end{array} \label{thick-plate}
\end{equation}
Again, the structure is the same as before with $\bu=(w,\btheta)$, the (generalized) strain $\bD\bu=(\boldsymbol{\chi},\boldsymbol{\gamma})$ consisting of the curvature $\boldsymbol{\chi}=\nabla \btheta$ and the shear strain $\boldsymbol{\gamma}=\nabla w -\btheta$. In the above, the support functions $\pi_G$ and $\Pi_G$ are defined with respect to a thick plate strength criterion involving both the bending moment tensor $\boldsymbol{M}$ and the shear force vector $\boldsymbol{Q}$ i.e. $\bSig=(\boldsymbol{M},\boldsymbol{Q})$. In \cite{bleyer2014stat_thick_plates,bleyer2015kin_thick_plates}, different choices of thick plate criteria are discussed, especially regarding the bending/shear interaction. For the sake of simplicity, we consider here only the case of no interaction between bending and shear, meaning that the thick plate strength criterion is in fact decoupled between bending and shear, taking the following form:
\begin{equation}
(\boldsymbol{M},\boldsymbol{Q})\in G_\text{thick plate} \Longleftrightarrow \begin{cases} \boldsymbol{M}\in G_\text{bend} \\ \boldsymbol{Q}\in G_\text{shear} \end{cases} \label{thick-plate-strength}
\end{equation}
where typically $G_\text{shear} =\{\boldsymbol{Q}\text{ s.t. } \|\boldsymbol{Q}\|_2\leq Q_0\}$ with $Q_0$ being the plate pure shear strength. Since typical thin plate criteria are SOC-representable, so will be the thick plate criterion \eqref{thick-plate-strength}. We also have the following expression for the corresponding support function:
\begin{align}
\pi_{G_\text{thick plate}}\left((\boldsymbol{\chi},\boldsymbol{\gamma})\right) &= \sup_{(\boldsymbol{M},\boldsymbol{Q})\in G_\text{thick plate}} \{\boldsymbol{M}:\boldsymbol{\chi}+\boldsymbol{Q}\cdot\boldsymbol{\gamma}\} \notag\\
&= \sup_{\boldsymbol{M}\in G_\text{bend}} \{\boldsymbol{M}:\boldsymbol{\chi}\} + \sup_{\boldsymbol{Q}\in G_\text{shear}} \{\boldsymbol{Q} \cdot\boldsymbol{\gamma}\} \notag\\
&= \pi_{G_\text{bend}}(\boldsymbol{\chi})+\pi_{G_\text{shear}}(\boldsymbol{\gamma})
\end{align}
with $\pi_{G_\text{shear}}(\boldsymbol{\gamma})=Q_0\|\boldsymbol{\gamma}\|_2$ and $\pi_{G_\text{bend}}$ depending on the choice of the bending strength criterion. Finally, the (generalized) velocity jump consisting of the rotation and velocity jumps $\bJ\bu=(\jump{\btheta}\overset{s}{\otimes}\bn,\jump{w}\bn)$ so that we also have:
\begin{equation}
\Pi_G\left((\jump{\btheta},\jump{w});\bn\right)=\pi_G\left(\left(\jump{\btheta}\overset{s}{\otimes}\bn,\jump{w}\bn\right)\right)
\end{equation}

\begin{figure}
\begin{center}
\includegraphics[width=0.6\textwidth]{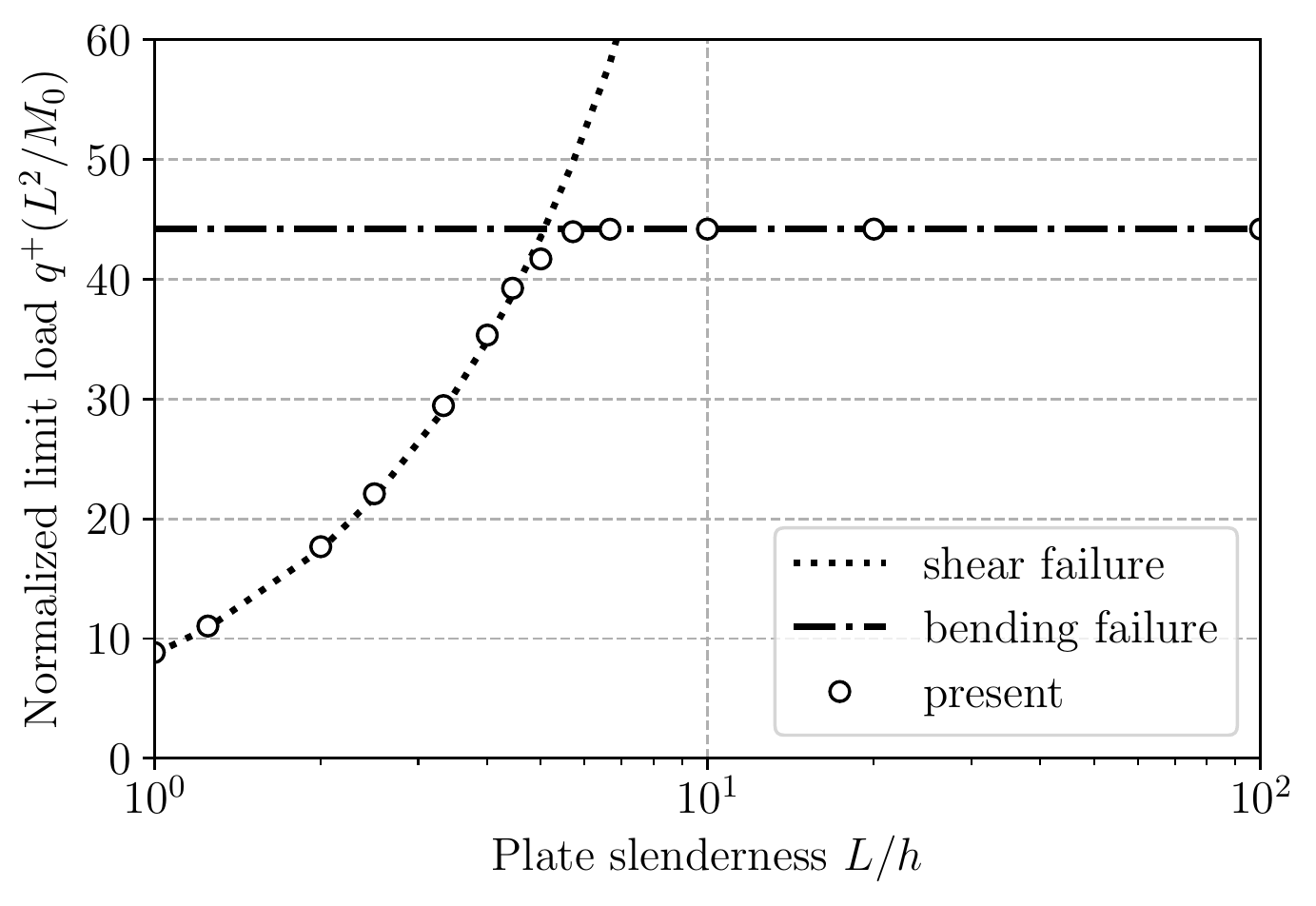}
\end{center}
\caption{Normalized limit load $q^+L^2/M_0$ of a square clamped plate as function of plate slenderness $L/h$: comparison with pure shear and pure bending failure mechanisms}
\label{thick-plate-res}
\end{figure}

\begin{figure}
\begin{center}
\begin{subfigure}{0.49\textwidth}
\includegraphics[width=\textwidth]{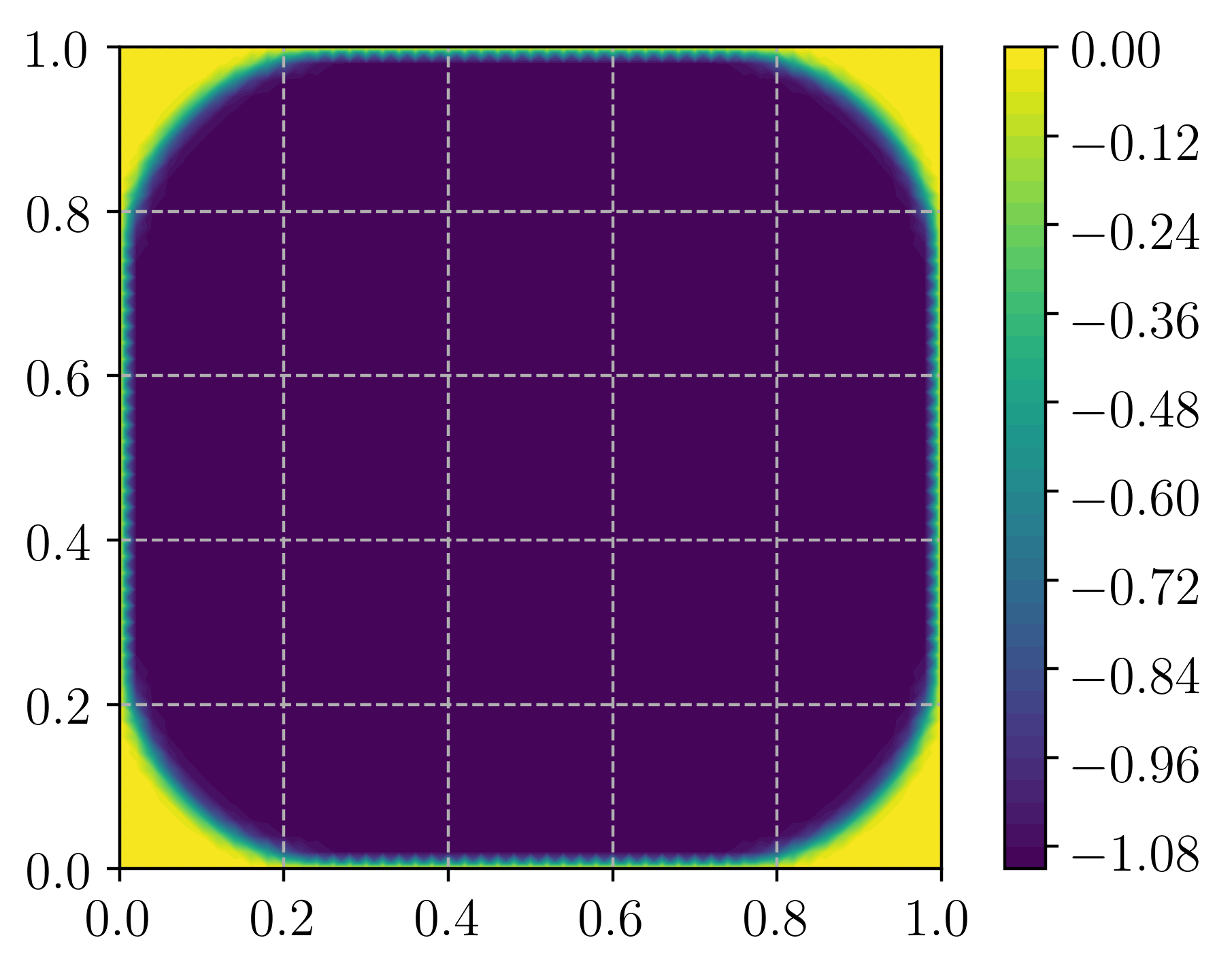}
\caption{Collapse mechanism}
\end{subfigure}
\hfill
\begin{subfigure}{0.49\textwidth}
\includegraphics[width=0.92\textwidth]{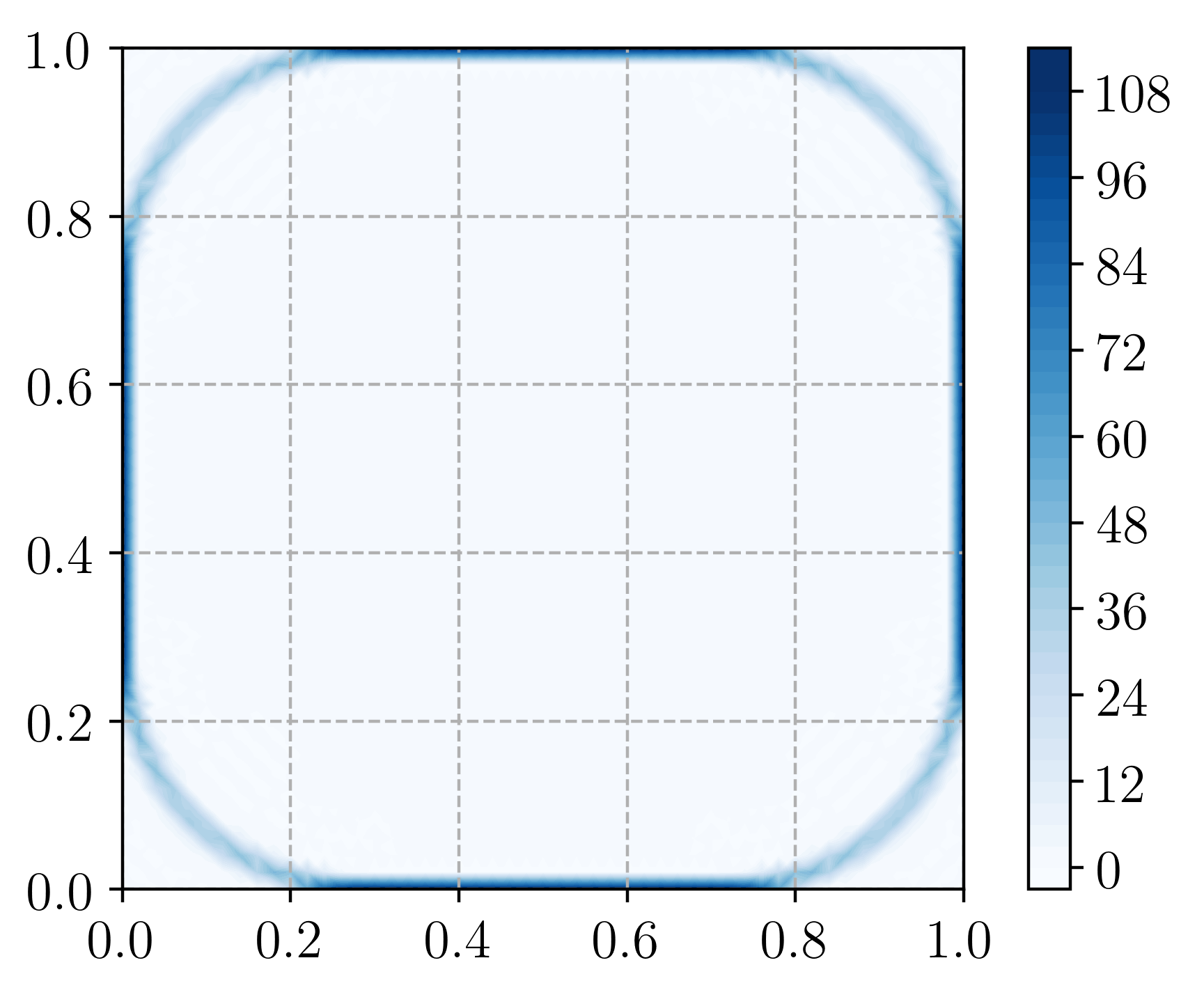}
\caption{Norm of the shear strain $\|\boldsymbol{\gamma}\|$}
\end{subfigure}
\end{center}
\caption{Pure shear collapse solution for a thick plate ($h/L=0.5$)}
\label{thick-shear}
\end{figure}

\begin{figure}
\begin{center}
\begin{subfigure}{0.49\textwidth}
\includegraphics[width=\textwidth]{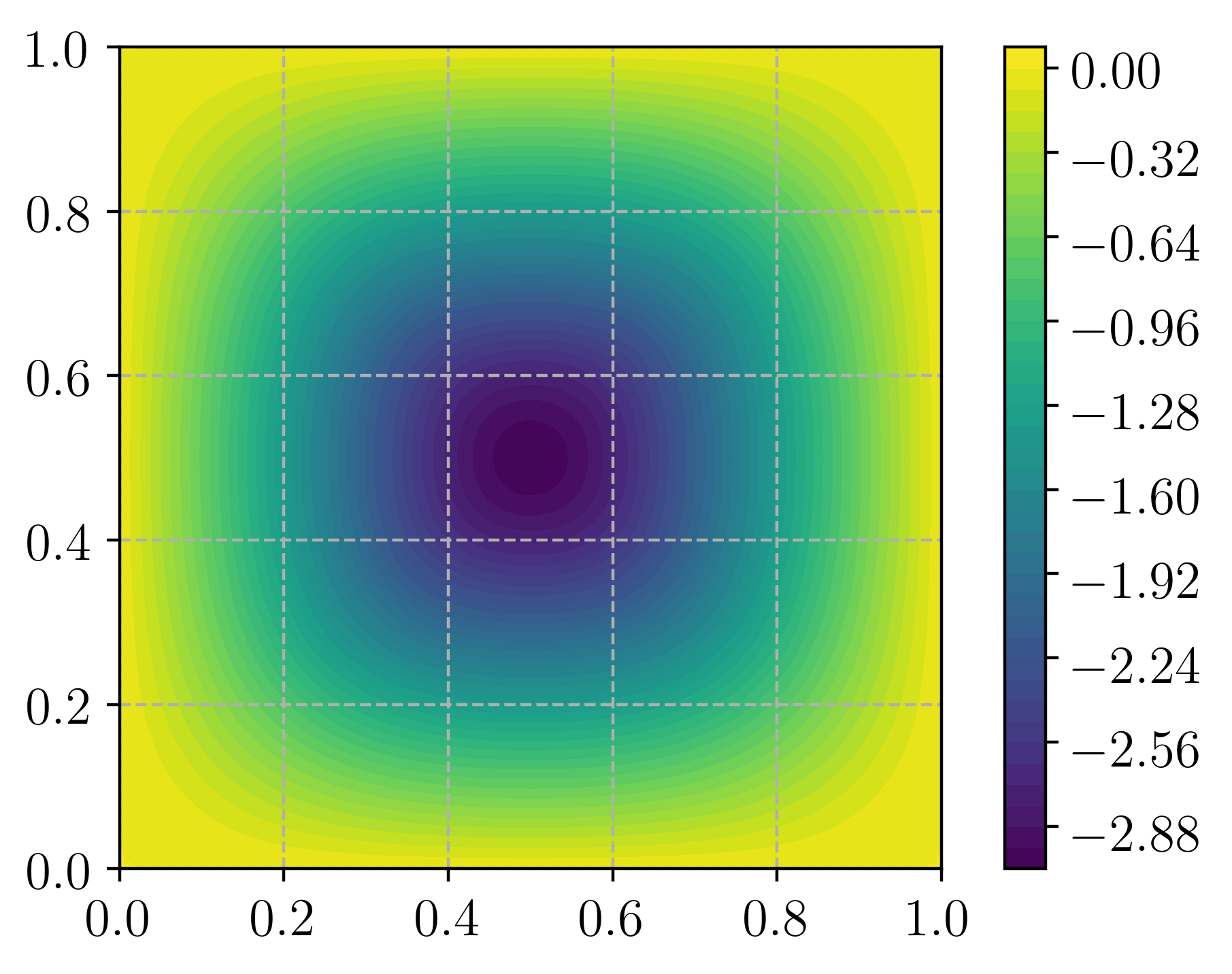}
\caption{Collapse mechanism}
\end{subfigure}
\hfill
\begin{subfigure}{0.49\textwidth}
\includegraphics[width=0.92\textwidth]{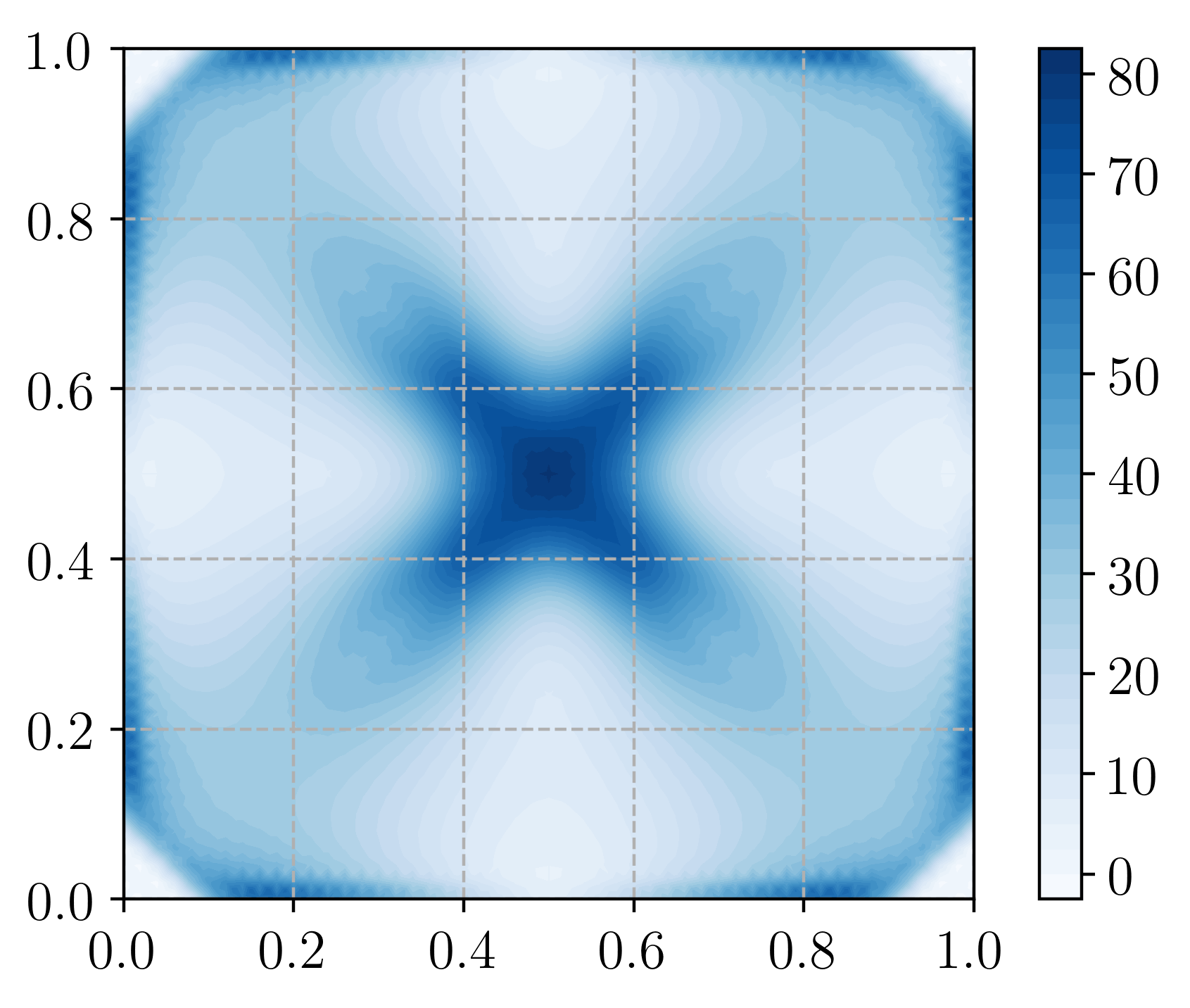}
\caption{Norm of the curvature $\|\boldsymbol{\chi}\|$}
\end{subfigure}
\end{center}
\caption{Pure bending collapse solution for a thin plate ($h/L=0.01$)}
\label{thick-bending}
\end{figure}

As discussed in \cite{bleyer2015kin_thick_plates}, finite-element discretization for the upper bound limit analysis of thick plates must be chosen with care. Indeed, continuous Lagrange interpolations for both the deflection $w$ and the rotation field $\btheta$ will exhibit shear locking in the thin plate limit. Reference \cite{bleyer2015kin_thick_plates} considered fully discontinuous interpolation of both fields and showed extremely good performances. To illustrate the versatility of the proposed framework, we consider here a small variant, namely a continuous $\PP^2$-Lagrange interpolation for the deflection $w$ and a discontinuous $\PP^1_d$-Lagrange interpolation for the rotation $\btheta$. This interpolation choice is expected to yield similar performances, although slightly higher, as the fully discontinuous $\PP^2_d/\PP^1_d$ interpolation of \cite{bleyer2015kin_thick_plates}. The problem of a uniformly distributed (intensity $q$) clamped square plate of side length $L$ and thickness $h$ is considered with a von Mises bending strength criterion (bending strength $M_0=\sigma_0 h^2/4$). The shear strength is taken as $Q_0=\sigma_0h/\sqrt{3}$. As already observed in \cite{bleyer2014stat_thick_plates,bleyer2015kin_thick_plates}, for the case of no bending/shear interaction, the limit load for this problem is well approximated by:
\begin{equation}
q^+ = \min\{q_\text{shear}^+;q_\text{bending}^+\}
\end{equation}
 where $q_\text{shear}^+=(2+\sqrt{\pi})Q_0/L$ is the pure shear solution for a square plate (see the Cheeger set example of \cite{bleyer2019automating}) and $q_\text{bending}^+\approx 44.2 M_0/L^2$ is the pure bending thin plate solution. This is indeed what is also observed for the present interpolation with a 50$\times$50 mesh, see Figure \ref{thick-plate-res}. The correct pure shear (Figure \ref{thick-shear}) and pure bending (Figure \ref{thick-bending}) mechanisms are also retrieved.


\subsection{Shells and multilayered plates}
Let us now briefly consider the case of shells, for which we will discuss here a lower bound implementation only. We refer to \cite{bleyer2016numerical} and references therein for more details concerning limit analysis of shells, especially the upper bound kinematic approach. The shell geometry will be approximated by a plane facet discretization into triangular elements and will be described locally by a unit normal $\udl{\nu}$ and a tangent plane spanned by two unit vectors $\udl{a}_1$ and $\udl{a}_2$. This local frame is therefore constant element-wise. The generalized internal forces of the shell model are described by a symmetric \textit{membrane force} tensor $\uddl{N}=N_{ij}\udl{a}_i\otimes\udl{a}_j$, a symmetric \textit{bending moment} tensor $\uddl{M}=M_{ij}\udl{a}_i\otimes\udl{a}_j$ and a \textit{shear force} vector $\udl{Q}=Q_i\udl{a}_i$ ($i,j=1,2$), the components of which are expressed in the local tangent plane. We will consider only thin shells, meaning that the shell strength criterion $G_\text{shell}$ is a convex set in the 6-dimensional $(\uddl{N},\uddl{M})$ space (infinite shear strength assumption). Introducing $\boldsymbol{T} = \uddl{N} + \udl{\nu}\otimes \udl{Q}$, the local equilibrium equations in a plane facet are given by:
\begin{align}
\div_T \boldsymbol{T} + \lambda\boldsymbol{f} &= 0 \label{shell-force-eq}\\
\div_T \uddl{M}+\udl{Q} &= 0 \label{shell-moment-eq}
\end{align}
where $\div_T$ is the tangent plane divergence operator and $\lambda \boldsymbol{f}$ a distributed loading with an amplification factor $\lambda$. In addition to local equilibrium, continuity equations of the force resultant $\boldsymbol{R}=\boldsymbol{T}\cdot\bn$ and the normal bending moment $\boldsymbol{\mathcal{M}}=\uddl{M}\cdot\udl{n}\times\udl{\nu}$ must be satisfied where $\bn$ is the in-plane normal to a facet edge:
\begin{align}
\jump{\boldsymbol{R}} &= 0 \label{shell-continuity-force}\\
\jump{\boldsymbol{\mathcal{M}}} &=0 \label{shell-continuity-moment}
\end{align}

Finite-element discretization relies on a discontinuous $\PP^1_d$ interpolation for both membrane forces $\uddl{N}$ and shear forces $\udl{Q}$ and a discontinuous $\PP^2_d$ interpolation for the bending moments $\uddl{M}$. Considering cell-wise uniform distributed loadings, local force equilibrium \eqref{shell-force-eq} is uniform and can therefore be satisfied exactly using a cell-wise constant Lagrange multiplier $\bu$. Local moment equilibrium \eqref{shell-moment-eq} can be satisfied exactly using a cell-wise linear Lagrange multiplier $\btheta$. Force and moment continuity equations \eqref{shell-continuity-force}-\eqref{shell-continuity-moment} are respectively satisfied by considering facet-wise linear and quadratic Lagrange multipliers $\boldsymbol{v}$ and $\boldsymbol{\vartheta}$. We again give below the main part of the corresponding Python code for formulating this rather complex problem:

\begin{pythoncode}
prob = MosekProblem("Shell lower bound limit analysis")
R = FunctionSpace(mesh, "R", 0)
Ne = VectorElement("DG", mesh.ufl_cell(), 1, dim=3)
Me = VectorElement("DG", mesh.ufl_cell(), 2, dim=3)
Qe = VectorElement("DG", mesh.ufl_cell(), 1, dim=2)
W = FunctionSpace(mesh, MixedElement([Ne, Me, Qe]))

lamb, Sig = prob.add_var([R, W])
prob.add_obj_func([1, None])

(N, M, Q) = split(Sig)
T = as_matrix([[N[0], N[2]],
               [N[2], N[1]],
               [Q[0], Q[1]]])
M = to_mat(M)
\end{pythoncode}
where the different unknown fields $\bSig=(\uddl{N},\uddl{M},\udl{Q})$ have been defined and collected into a global vector. Equilibrium and continuity equations are then defined as in section \ref{sec:static-continuum} for continua by specifying the Lagrange multiplier function space and writing the weak form of the constraint:
\begin{pythoncode}
V_f_eq = VectorFunctionSpace(mesh, "DG", 0, dim=3)
def force_equilibrium(u):
    u_loc = dot(Ploc, u)
    return [dot(u, f)*lamb*dx, dot(u_loc, divT(T))*dx]
prob.add_eq_constraint(V_f_eq, A=force_equilibrium)

V_m_eq = VectorFunctionSpace(mesh, "DG", 1, dim=2)
def moment_equilibrium(theta):
    return [None, dot(theta, divT(M)+Q)*dx]
prob.add_eq_constraint(V_m_eq, A=moment_equilibrium)

V_f_jump = VectorFunctionSpace(mesh, "Discontinuous Lagrange Trace", 1, dim=3)
Tglob = dot(Ploc.T, T)
def force_continuity(v):
    return [None, dot(avg(v),jump(Tglob, n_plan))*dS]
prob.add_eq_constraint(V_f_jump, A=force_continuity)

V_m_jump = VectorFunctionSpace(mesh, "Discontinuous Lagrange Trace", 2, dim=3)
Mglob = dot(Ploc_plane.T, M)
def moment_continuity(vtheta):
    return [None, dot(avg(vtheta), cross(jump(Mglob, n_plan), avg(nu)))*dS]
prob.add_eq_constraint(V_m_jump, A=moment_continuity)
\end{pythoncode}
where the \texttt{Ploc} (resp. \texttt{Ploc\_plane}) variable corresponds to a rotation matrix transforming fields expressed in the global $(\boldsymbol{e}_x,\boldsymbol{e}_y,\boldsymbol{e}_z)$ into the local $(\udl{a}_1,\udl{a}_2,\udl{\nu})$ (resp. $(\udl{a}_1,\udl{a}_2)$) frame. Note that we approximated the bending moment continuity equation by $\jump{\uddl{M}\cdot\udl{n}}\times\widehat{\udl{\nu}}$ using an average normal vector $\widehat{\udl{\nu}} = (\udl{\nu}^++\udl{\nu}^-)/2$.

If the previous snippets illustrate the efficiency of FEniCS high-level symbolic formulations for this kind of complex problem, the conic representation framework will also be extremely beneficial for formulating the shell strength criterion. Indeed, as discussed in length in \cite{bleyer2016numerical}, even for the simple case of a homogeneous von Mises thin shell, the strength condition expressed in terms of $(\uddl{N},\uddl{M})$ stress-resultant becomes extremely complicated \cite{ilyushin1956plasticite} so that simple SOC-representable approximations have been proposed in the past for the von Mises shell \cite{robinson1971comparison}. In \cite{bleyer2016numerical}, we proposed a general way of formulating an $(\uddl{N},\uddl{M})$  shell strength criterion for a general multilayered shell through an up-scaling procedure. It is given by:
\begin{equation}
(\uddl{N},\uddl{M})\in G_\text{shell} \Longleftrightarrow \begin{cases}
\exists\: \bsig(z) \in \mathcal{G}_\text{ps} \quad \forall z\in[-h/2;h/2] \text{ and s.t.} \\
\displaystyle{\uddl{N} = \int_{-h/2}^{h/2}\bsig(z) \dz} \\
\displaystyle{\uddl{M} = \int_{-h/2}^{h/2}(-z)\bsig(z) \dz}
\end{cases} \label{G-shell}
\end{equation}
where $h$ is the shell thickness and $\mathcal{G}_{ps}$ is the material local plane stress criterion, which may potentially depend on coordinate $z$ for a multilayered shell. To make formulation \eqref{G-shell} usable in practice, the local plane stress distribution $\bsig(z)$ is replaced by a discrete set of plane stress states $\bsig_g =\bsig(z_g)$ expressed at quadrature points $z_g$ which are used to approximate the two integrals:
\begin{equation}
(\uddl{N},\uddl{M})\in G_\text{shell}^\text{approx} \Longleftrightarrow \begin{cases}
\exists\: \bsig_g \in \mathcal{G}_\text{ps} \quad \forall g=1,\ldots,n_z \text{ and s.t.} \\
\displaystyle{\uddl{N} = \sum_{g=1}^{n_z}\omega_g \bsig_g} \\
\displaystyle{\uddl{M} = \sum_{g=1}^{n_z}(-z_g)\omega_g \bsig_g}
\end{cases} \label{G-shell-approx}
\end{equation}
where $\omega_g$ are the corresponding quadrature weights of the $n_z$-points quadrature rule. As discussed in \cite{bleyer2016numerical}, the precise choice of the quadrature rule leads to different kinds of approximations to $G_\text{shell}$: e.g. an upper bound approximation is obtained with a trapezoidal quadrature rule, a Gauss-Legendre quadrature leads to an approximation with no lower or upper bound status, a rectangular rule will give a lower bound approximation. In the following, we choose the latter to be consistent with the lower bound status of the static approach.

From \eqref{G-shell-approx} it can readily be seen that if $\mathcal{G}_\text{ps}$ is SOC-representable, so will be $G_\text{shell}^\text{approx}$. For instance, in the case of a plane stress von Mises criterion of uniaxial strength $\sigma_0$, it can be shown that:
\begin{equation}
\bsig \in \mathcal{G}_\text{ps von Mises} \Leftrightarrow \begin{cases} \bsig = \mathbf{J}\overline{\mathbf{y}} \\
y_0 = \sigma_0 \\
\|\overline{\mathbf{y}}\|\leq y_0
\end{cases} \quad \text{where } \mathbf{J} = \begin{bmatrix}
1 & 1/\sqrt{3} & 0 \\
0 & 2/\sqrt{3} & 0 \\
0 & 0 & 1/\sqrt{3}
\end{bmatrix}
\end{equation}
in which the last constraint is a $\Qq_4$ quadratic cone. Choosing a rectangular quadrature rule with $\omega_g=h/n_z$ and $z_g = -\frac{h}2 +\frac{h}{n_z}(g-\frac{1}{2})$ for $g=1,\ldots, n_z$, we therefore have:
\begin{equation}
(\uddl{N},\uddl{M})\in G_\text{shell}^\text{approx} \Longleftrightarrow \begin{cases}
\exists\: \mathbf{y}_g=(y_{0g},\overline{\mathbf{y}}_g) \in \Qq_4 \quad \forall g=1,\ldots,n_z \text{ and s.t.} \\
y_{0g} = \sigma_0 \\
\displaystyle{\uddl{N} = \sum_{g=1}^{n_z}\frac{h}{n_z} \mathbf{J}\overline{\mathbf{y}}_g} \\
\displaystyle{\uddl{M} = \sum_{g=1}^{n_z}(-z_g)\frac{h}{n_z} \mathbf{J}\overline{\mathbf{y}}_g}
\end{cases} \label{G-vonMises-approx}
\end{equation}
which obviously fits format \eqref{conic-representable}. Let us finally remark that the approximation will converge to the shell criterion $G_\text{shell}$ when increasing $n_z$. In the following we took $n_z=6$.

\begin{figure}
\begin{center}
\begin{subfigure}{0.54\textwidth}
\includegraphics[width=\textwidth]{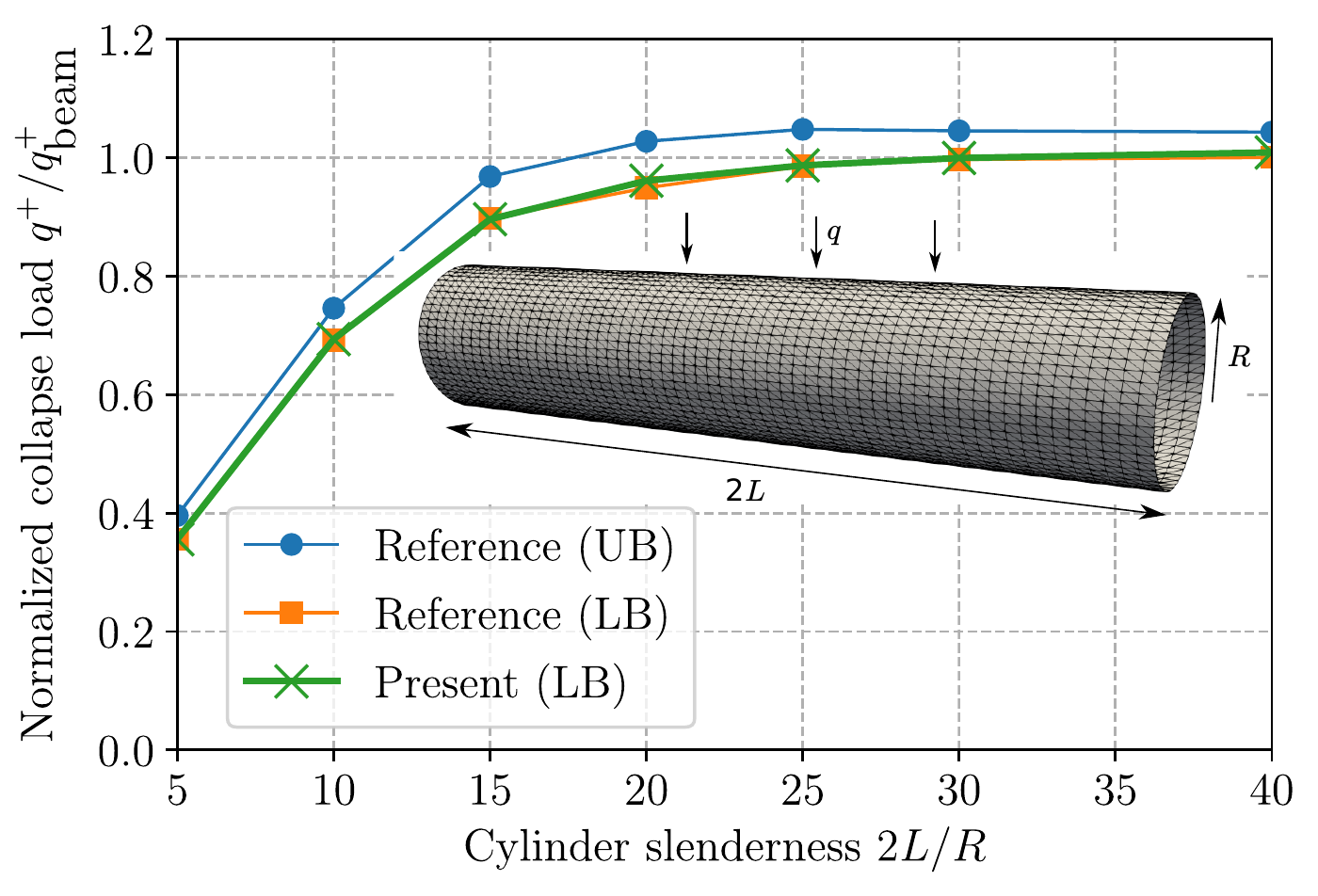}
\caption{Limit load results: reference LB and UB limit loads from \cite{bleyer2016numerical}}
\label{shell-results}
\end{subfigure}
\hfill
\begin{subfigure}{0.45\textwidth}
\includegraphics[width=\textwidth]{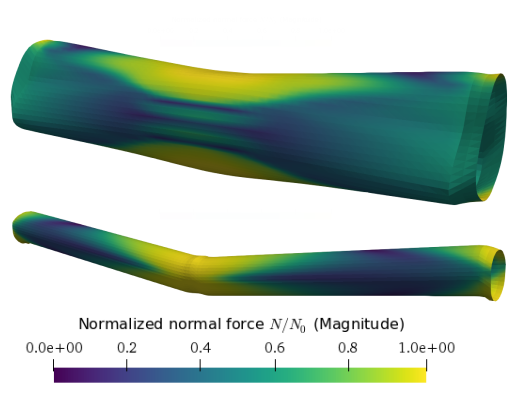}
\caption{Cylindrical shell pseudo-collapse mechanisms and normalized normal force magnitude $\|\boldsymbol{N}\|/N_0$. Top: slenderness $2L/R=10$, bottom: slenderness $2L/R=30$.}
\label{shell-mechanisms}
\end{subfigure}
\end{center}
\caption{Cylindrical shell under self-weight}
\end{figure}

As an illustrative application, we consider the problem of a cylindrical shell of length $2L$, radius $R$ and thickness $h=0.01R$, clamped at both extremities and loaded by a self-weight uniform vertical loading $\boldsymbol{f}=-q\boldsymbol{e}_z$ (see Figure \ref{shell-results}). The shape of the collapse mechanism varies depending on the cylinder slenderness $2L/R$. For sufficiently long cylinders, the computed limit load $q^+$ is well described by the one obtained when representing the cylinder as a 1D beam $q_\text{beam}^+ = \frac{32}{\pi}N_0\left(\dfrac{R}{2L}\right)^2$ with $N_0=\sigma_0 h$ being the membrane uniaxial strength. The obtained limit loads agree very well with the lower bound result of \cite{bleyer2016numerical}. Having access to the Lagrange multiplier $\bu$, we reconstructed a pseudo-collapse mechanism by performing a projection of $\bu$ on a continuous $\PP^1$ space. The obtained deformations have been represented in Figure \ref{shell-mechanisms} along with the normal force magnitude distribution. It can be seen that the mechanisms agree well with those obtained from an upper bound kinematic approach in \cite{bleyer2016numerical}, with a beam-like mechanism involving plastic hinges at the clamped supports and mid-span for the case $2L/R=30$.

\section{\label{sec:gen-cont}Generalized continua}
In this last section, we further illustrate the proposed framework on two generalized continuum models, namely a strain gradient and a Cosserat continuum. We would like to point out that the numerical implementation of limit analysis for these kinds of model is almost non-existent whereas we will show that it can now be easily formulated with the proposed framework.

\subsection{Strain-gradient material}
In this subsection, we consider the extension of limit analysis to strain gradient materials. We do not attempt at providing physical justifications for using this kind of model but let us just mention that it can, for instance, be obtained when considering the elastically rigid version of a strain-gradient plasticity model. We therefore consider the following strain-gradient generalization of the kinematic limit analysis theorem:
 \begin{equation}
\begin{array}{rl}
\displaystyle{\inf_{\boldsymbol{u}\in \Vv_h}} & \displaystyle{\int_{\Omega} \pi_{G_\text{SG}}((\nabla^s \boldsymbol u, \nabla^2 \bu)) \dx} \label{kinematic-strain-gradient}\\
\text{s.t.} & \int_{\Omega}\boldsymbol{f}\cdot\boldsymbol{u}\dx + \int_{\partial \Omega_T} \boldsymbol{t}\cdot\boldsymbol{u}\dS = 1
\end{array}
\end{equation}
where we considered only classical loadings (body forces or surface tractions) and where the strain operator $\bD\bu = (\nabla^s \bu ,\nabla^2\bu)$ now includes both the first and second displacement gradient of $\bu$ with $(\nabla^2\bu)_{ijk}=u_{i,jk}$. We do not consider here the equivalent static formulation but let us just point out that the generalized stress $\bSig=(\bsig, \boldsymbol{\tau})$ includes both the standard Cauchy stress $\bsig$ and the third-rank couple stress tensor $\boldsymbol{\tau}=(\tau_{ijk})$ which is associated by duality with the second gradient $\nabla^2 \bu$. The generalized strength criterion therefore depends both on $\bsig$ and $\boldsymbol{\tau}$. For simplicity, we will consider the following extended von Mises criterion:
\begin{equation}
(\bsig,\boldsymbol{\tau})\in G_\text{SG} \Leftrightarrow \sqrt{\frac{1}{2}(\boldsymbol{s}:\boldsymbol{s}+\ell^{-2}\boldsymbol{\tau}\T\tripleC\boldsymbol{\tau})} \leq k
\end{equation}
where $\boldsymbol{s}=\dev\bsig$, $\boldsymbol{\tau}\T\tripleC\boldsymbol{\tau}=\tau_{ijk}\tau_{ijk}$ and $\ell$ is an internal length scale. The associated support function is:
\begin{equation}
\pi_{G_\text{SG}}((\boldsymbol{d},\boldsymbol{\eta})) = \sup_{(\bsig,\boldsymbol{\tau})\in G_\text{SG}}\{\bsig:\boldsymbol{d}+\boldsymbol{\tau}\T\tripleC\boldsymbol{\eta} \}=\begin{cases} k\sqrt{2(\boldsymbol{d}:\boldsymbol{d}+\ell^2\boldsymbol{\eta}\T\tripleC\boldsymbol{\eta})} & \text{if } \tr\boldsymbol{d}=0 \\
+\infty & \text{otherwise} \end{cases} \label{pi-vM-SG}
\end{equation}
where $\boldsymbol{d}=\nabla^s \bu$ and $\boldsymbol{\eta}=\nabla^2\bu$.

Restricting to a plane strain situation, one has $\eta_{ij3}=\eta_{i3j}=\eta_{3ij}=0$ and:
\begin{align}
\eta_{111}=u_{1,11} &\qquad \eta_{211}=u_{2,11} \notag\\
\eta_{122}=u_{1,22} &\qquad   \eta_{222}=u_{2,22}\\
\eta_{112}=\eta_{121}=u_{1,12}  &\qquad \eta_{212}=\eta_{221}=u_{2,12} \notag
\end{align}
Introducing $\mathbf{D}=(d_{11},d_{22},\sqrt{2}d_{12},\ell \eta_{111},\ell\eta_{122}, \sqrt{2}\ell \eta_{112}, \ell\eta_{211}, \ell\eta_{222}, \sqrt{2}\ell\eta_{212})$, one has $\mathbf{D}\T\mathbf{D} = \boldsymbol{d}:\boldsymbol{d}+\ell^2\boldsymbol{\eta}\T\tripleC\boldsymbol{\eta}$ so that $\pi_{G_\text{SG}}((\boldsymbol{d},\boldsymbol{\eta})) = k\sqrt{2}\|\mathbf{D}\|_2$. Since $\pi_{G_\text{SG}}$ involves a $L_2$-norm on a 9-dimensional vector, it can be represented using a 10-dimensional quadratic cone $\Qq_{10}$ (see \cite{bleyer2019automating}).

As regards the finite-element discretization, we choose a $\PP^2$-Lagrange interpolation for $\bu$. In formulation \eqref{kinematic-strain-gradient}, it is implicitly assumed that both $\bu$ and $\nabla \bu$ are continuous. The latter condition will not be achieved easily by a standard FE discretization so that we supplement \eqref{kinematic-strain-gradient} by a discontinuity term for $\partial_n \bu = \nabla \bu \cdot \bn$, similarly to thin plates:
 \begin{equation}
\begin{array}{rl}
\displaystyle{\inf_{\boldsymbol{u}\in \Vv_h}} & \displaystyle{\int_{\Omega} \pi_{G_\text{SG}}((\nabla^s \boldsymbol u, \nabla^2 \bu)) \dx + \int_{\Gamma} \pi_{G_\text{SG}}((0, \jump{\partial_n \bu} \otimes \bn\otimes\bn)) \dS} \label{kinematic-strain-gradient-disc}\\
\text{s.t.} & \int_{\Omega}\boldsymbol{f}\cdot\boldsymbol{u}\dx + \int_{\partial \Omega_T} \boldsymbol{t}\cdot\boldsymbol{u}\dS = 1
\end{array}
\end{equation}
with $\pi_{G_\text{SG}}((0, \partial_n \bu \otimes \bn\otimes\bn)) = k\sqrt{2}\ell\|\jump{\partial_n \bu}\|_2$ for \eqref{pi-vM-SG}. Again, this can be easily implemented in very few lines of code, regarding the problem complexity:
\begin{pythoncode}
prob = MosekProblem("Strain gradient limit analysis")
u = prob.add_var(V, bc=bc)

D = as_vector([u[0].dx(0), u[1].dx(1), (u[0].dx(1)+u[1].dx(0))/sqrt(2),
               l*u[0].dx(0).dx(0), l*u[0].dx(1).dx(1), sqrt(2)*l*u[0].dx(0).dx(1),
               l*u[1].dx(0).dx(0), l*u[1].dx(1).dx(1), sqrt(2)*l*u[1].dx(0).dx(1)])
pi = L2Norm(D, quadrature_scheme="vertex")
prob.add_convex_term(k*sqrt(2)*pi)

isochoric = EqualityConstraint(div(u), quadrature_scheme="vertex")
prob.add_convex_term(isochoric)

n = FacetNormal(mesh)
pi_d = L2Norm([jump(k*sqrt(2)*l*grad(u), n)], on_facet=True)
prob.add_convex_term(pi_d)

prob.optimize()
\end{pythoncode}

\begin{figure}
\begin{center}
\includegraphics[width=0.6\textwidth]{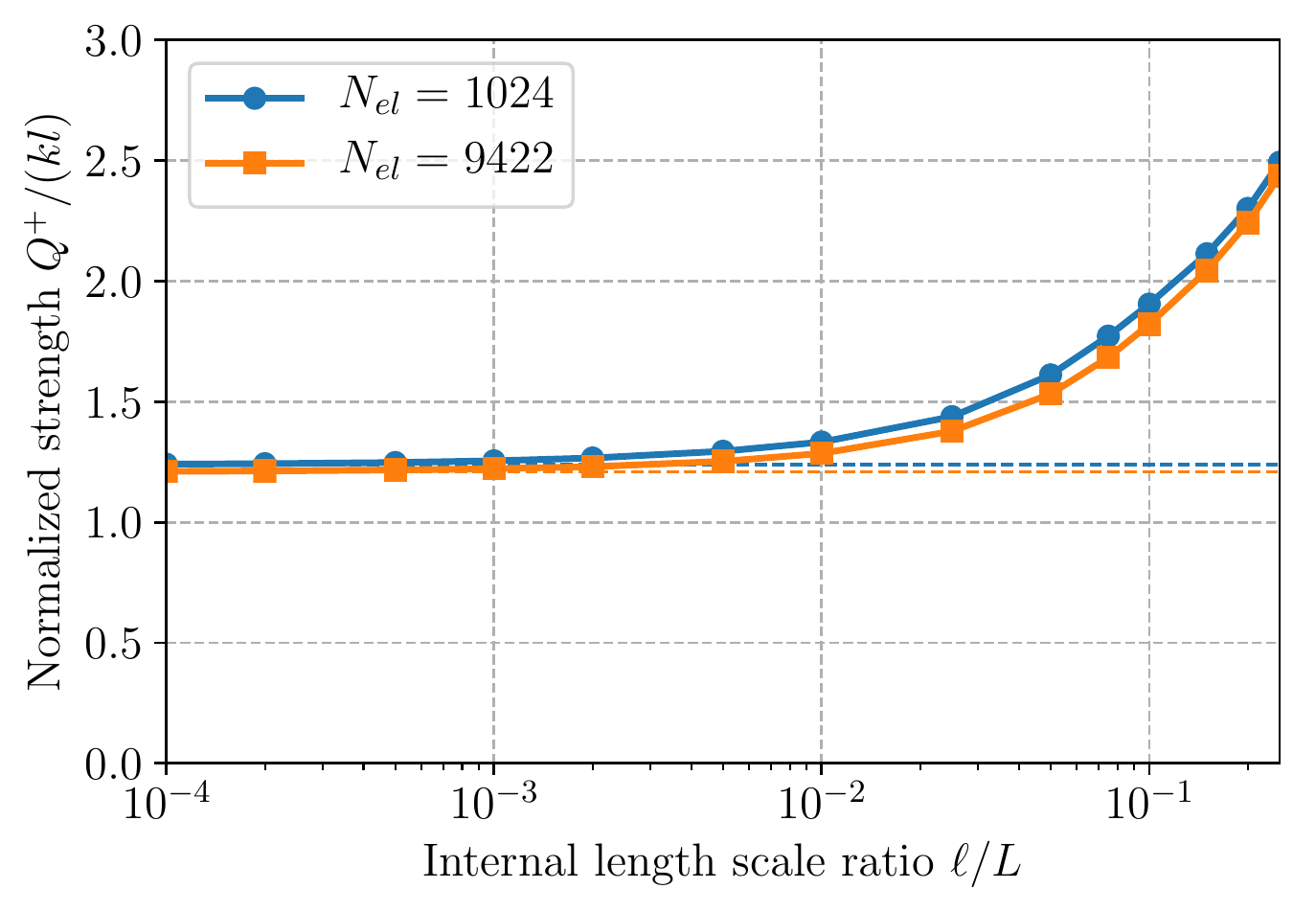}
\end{center}
\caption{Evolution of the normalized uniaxial strength for the strain-gradient perforated traction as function of the internal length scale $\ell/L$ for two mesh sizes (dashed lines correspond to standard continuum limit analysis results $\ell=0$).}
\label{strain-gradient-results}
\end{figure}

\begin{figure}
\begin{center}
\begin{subfigure}[t]{0.24\textwidth}
\includegraphics[width=\textwidth]{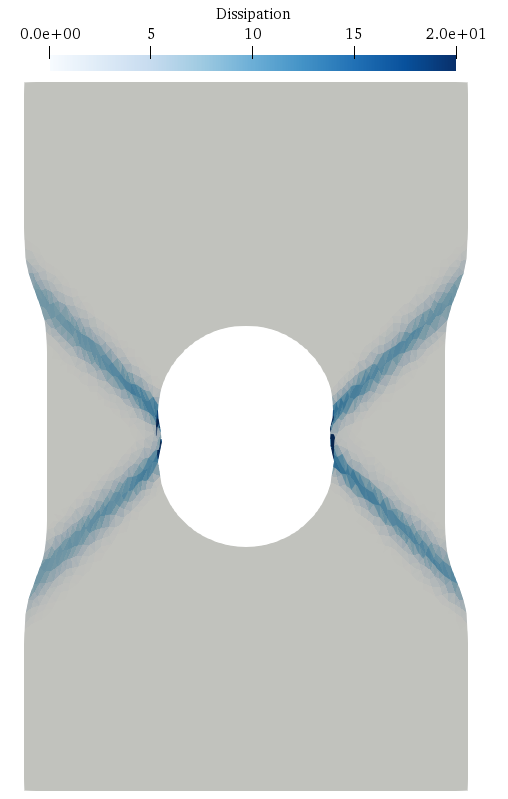}
\caption{$\ell/L=0$ (standard continuum)}
\end{subfigure}
\hfill
\begin{subfigure}[t]{0.24\textwidth}
\includegraphics[width=\textwidth]{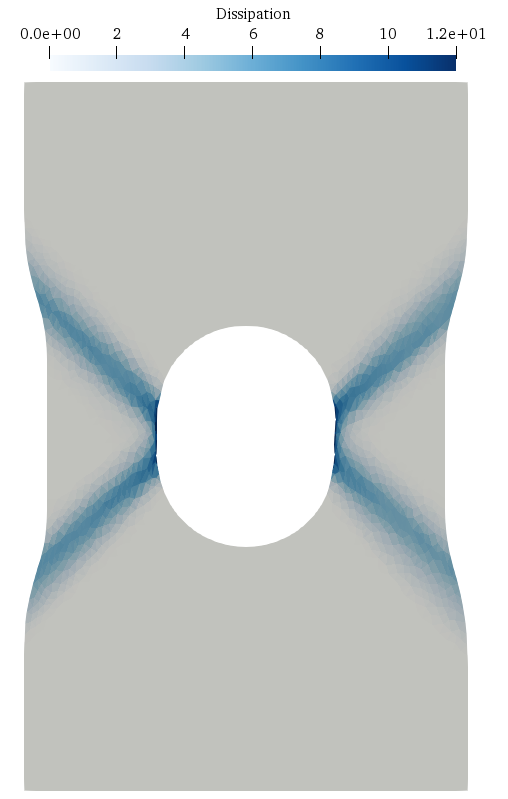}
\caption{$\ell/L=0.001$}
\end{subfigure}
\hfill
\begin{subfigure}[t]{0.24\textwidth}
\includegraphics[width=\textwidth]{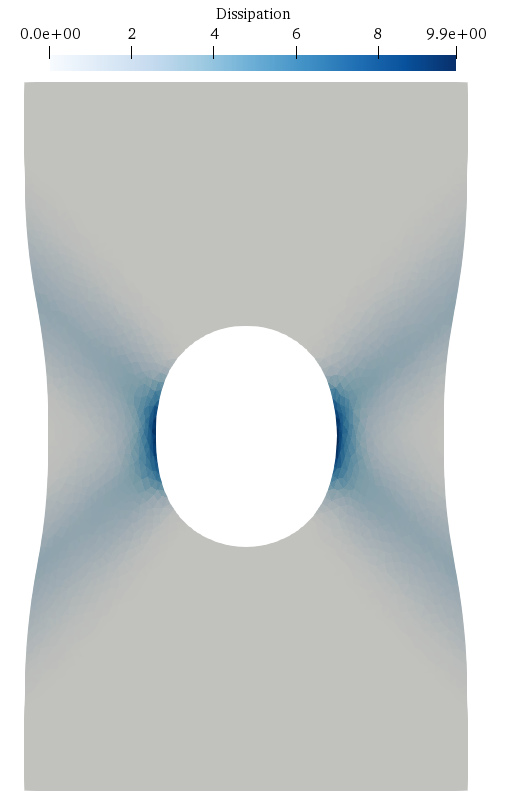}
\caption{$\ell/L=0.01$}
\end{subfigure}
\hfill
\begin{subfigure}[t]{0.24\textwidth}
\includegraphics[width=\textwidth]{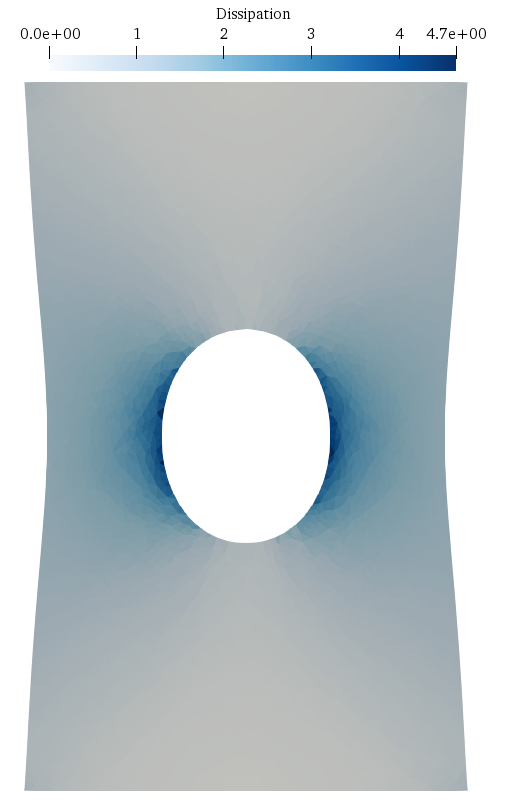}
\caption{$\ell/L=0.1$}
\end{subfigure}
\end{center}
\caption{Collapse mechanism and plastic dissipation as function of internal length scale ratio $\ell/L$}
\label{strain-gradient-mechanisms}
\end{figure}

As an illustrative application, we consider a rectangular domain of dimensions $L\times 1.5L$ perforated at its center by a circular hole of radius $R=0.2L$. The bottom boundary is fully clamped and the top one is displaced vertically $\bu=(0,U)$. No other loading is applied and the computed objective value of \eqref{kinematic-strain-gradient-disc} will be $Q^+U$ where $Q=\int_{y=H}\sigma_{yy}\dS$ is the resultant force on the top boundary. The evolution of the normalized uniaxial strength $Q^+/(kL)$ is plotted for various values of the internal length parameter $\ell/L$ in Figure \ref{strain-gradient-results} for two different mesh sizes. As expected, the plate apparent strength is size-dependent and exhibits a strengthening behaviour for larger values of $\ell$ or, equivalently, smaller sample size $L$. The standard continuum limit analysis results (dashed lines) are retrieved when $\ell/L\to 0$. Collapse mechanisms along with dissipation fields $\pi_\text{SG}$ are represented in Figure \ref{strain-gradient-mechanisms}. Broadening of the plastic dissipation slip zones can clearly be observed for increasing values of $\ell/L$.
  

\subsection{A Cosserat-continuum model for jointed rocks}

\begin{figure}
\begin{center}
\includegraphics[width=0.5\textwidth]{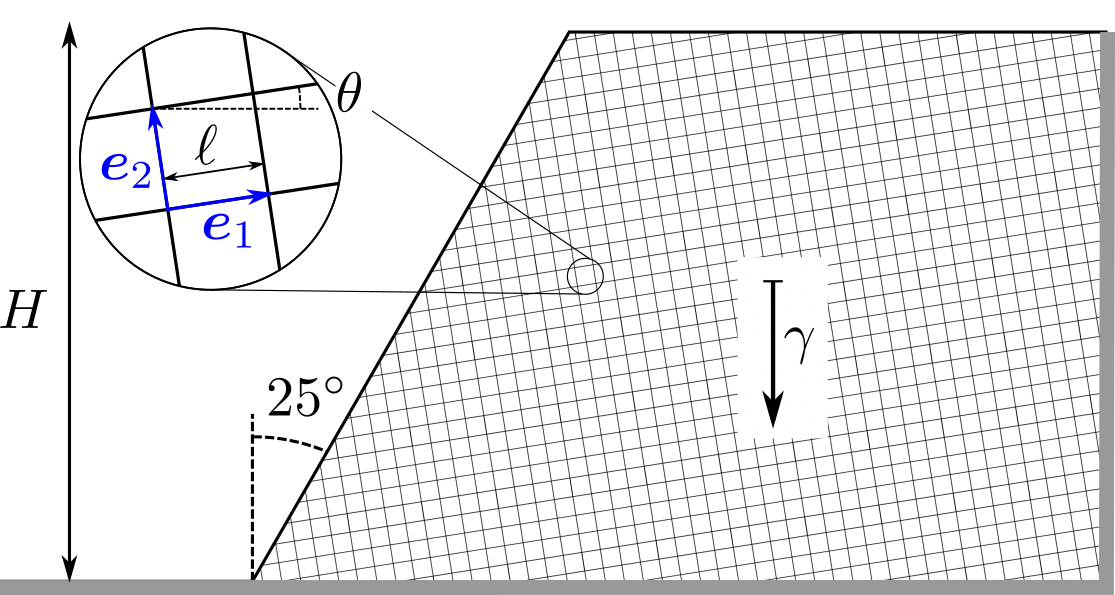}
\end{center}
\caption{Stability of a jointed rock excavation}
\label{jointed-rock-geo}
\end{figure}

We further illustrate the ability of the proposed framework to tackle generalized continua by considering a Cosserat (or micropolar continuum) model for jointed rocks, initially proposed in \cite{de2002failure}. The governing equations of the model, in plane strain conditions, involve a non-symmetric stress tensor $\bSig = \begin{bmatrix}
\Sigma_{11} & \Sigma_{12} \\ \Sigma_{21} & \Sigma_{22}
\end{bmatrix}$ and a couple stress vector $\boldsymbol{H} = (H_1, H_2)$ both expressed in the local reference frame $(\boldsymbol{e}_1,\boldsymbol{e}_2)$ of the jointed rock mass (see Figure \ref{jointed-rock-geo}). The corresponding equilibrium equations read as:
\begin{align}
\div \bSig + \boldsymbol{f} &= 0 \\
\div \boldsymbol{H}+\Sigma_{21}-\Sigma_{12} &=0 \label{Cosserat-eq2}
\end{align}
the corresponding weak form obtained from the virtual work principle being:
\begin{equation}
\int_\Omega \left(\bSig\T:(\nabla\bu -\skew \omega)+\boldsymbol{H}\cdot\nabla \omega \right)\dx = \int_{\Omega}\boldsymbol{f}\cdot\boldsymbol{u}\dx
\end{equation}
for any continuous test function $\bu$ and $\omega$ with $\skew\omega = \omega(\udl{e}_2\otimes\udl{e}_1-\udl{e}_1\otimes\udl{e}_2)$ and where we considered only body forces as loading parameters.

As regards strength properties, the rock mass is assumed to obey a Mohr-Coulomb criterion of cohesion $c_m$ and friction angle $\phi_m$. The joints are represented as an orthogonal array, spaced by a length $\ell$ and making an angle $\theta$ with the horizontal axis. They are assumed to also obey a Mohr-Coulomb condition with parameters $(c_j,\phi_j)$. The generalized strength condition for a jointed rock mass modelled as a Cosserat continuum is expressed as \cite{de2002failure}: 
\begin{equation}
(\bSig,\boldsymbol{H})\in G_\text{Cosserat} \Longleftrightarrow \begin{cases}
\Sigma_{11}\tan\phi_j+|\Sigma_{21}|\leq c_j \\
\Sigma_{22}\tan\phi_j+|\Sigma_{12}|\leq c_j \\
\dfrac{\ell}{2}\Sigma_{11}+|H_1|\leq \dfrac{\ell}{2}\dfrac{c_j}{\tan\phi_j} \\
\dfrac{\ell}{2}\Sigma_{22}+|H_2|\leq \dfrac{\ell}{2}\dfrac{c_j}{\tan\phi_j} \\
|H_1|\leq \dfrac{\ell}{2}\dfrac{c_j}{\cos\phi_j} \\
|H_2|\leq \dfrac{\ell}{2}\dfrac{c_j}{\cos\phi_j} \\
\sym\bSig \in G_\text{MC,2D}(c_m,\phi_m)
\end{cases} \label{joint-rock-G}
\end{equation}
where the last condition expresses the rock mass Mohr-Coulomb criterion on $\sym\bSig = (\bSig+\bSig\T)/2$ and where all other conditions involve the joints resistance. Let us point out that the case $\ell=0$ induces $H_i=0$ and thus $\bSig = \bSig\T$ due to \eqref{Cosserat-eq2}, one therefore retrieves a Cauchy model with a strength criterion described by the first, second and last conditions of \eqref{joint-rock-G}. Finally, $G_\text{Cosserat}$ involves only linear inequality constraints in addition to the Mohr-Coulomb criterion $G_\text{MC,2D}$. It is, therefore, SOC-representable, the part involving joints only being linear-representable.\\

\begin{figure}
\begin{center}
\includegraphics[width=0.6\textwidth]{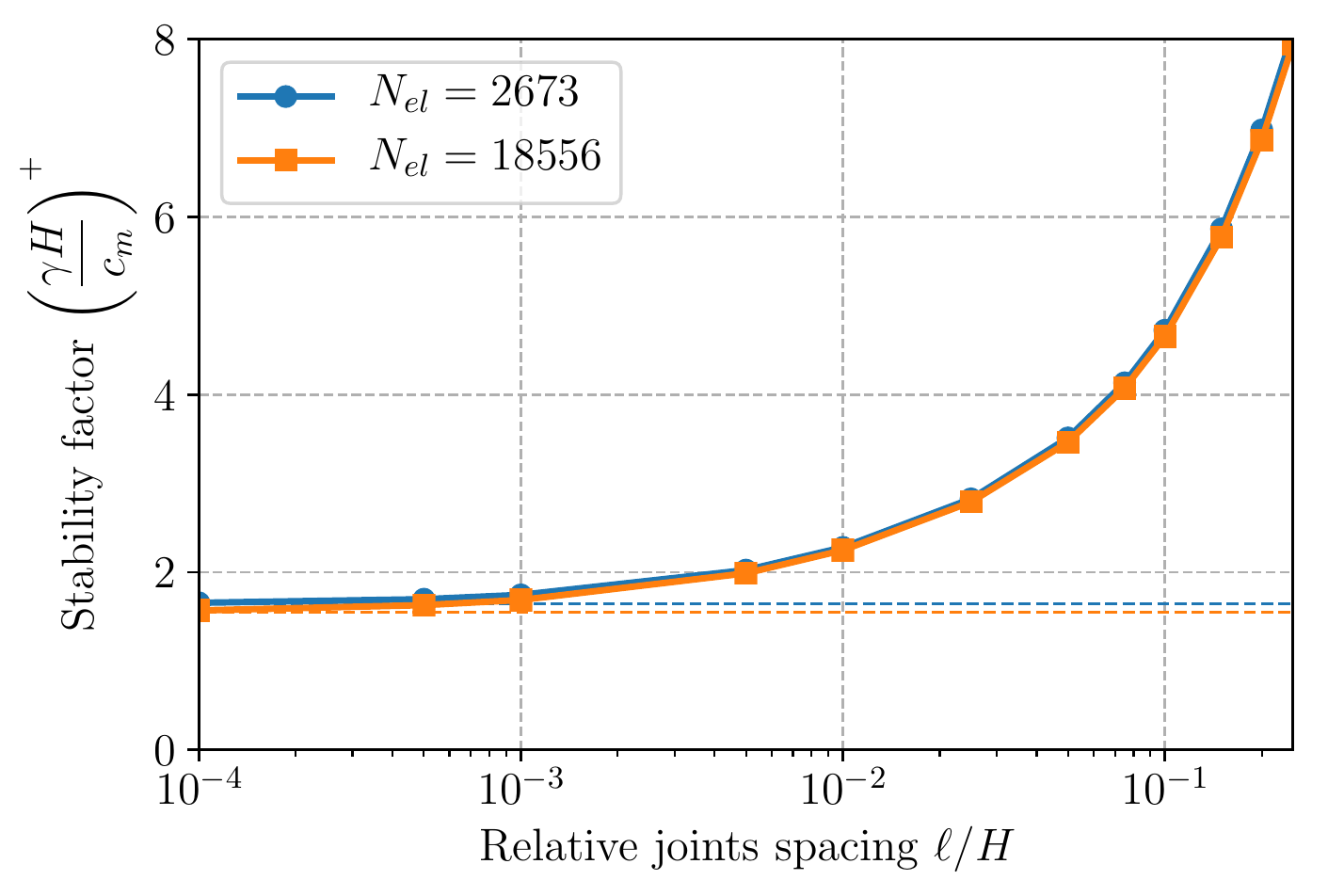}
\end{center}
\caption{Evolution of the stability factor as function of the joints spacing $\ell/H$ for two mesh sizes (dashed lines correspond to standard continuum limit analysis results $\ell=0$).}
\label{cosserat-results}
\end{figure}

\begin{figure}
\begin{center}
\begin{subfigure}{0.32\textwidth}
\includegraphics[width=\textwidth]{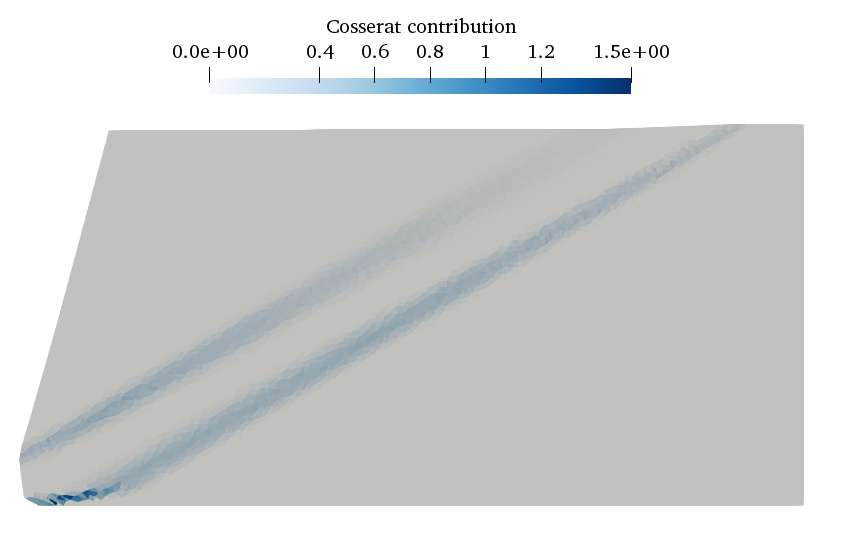}
\caption{$\ell/H=0.001$}
\label{cosserat-mech-1e-3}
\end{subfigure}
\hfill
\begin{subfigure}{0.32\textwidth}
\includegraphics[width=\textwidth]{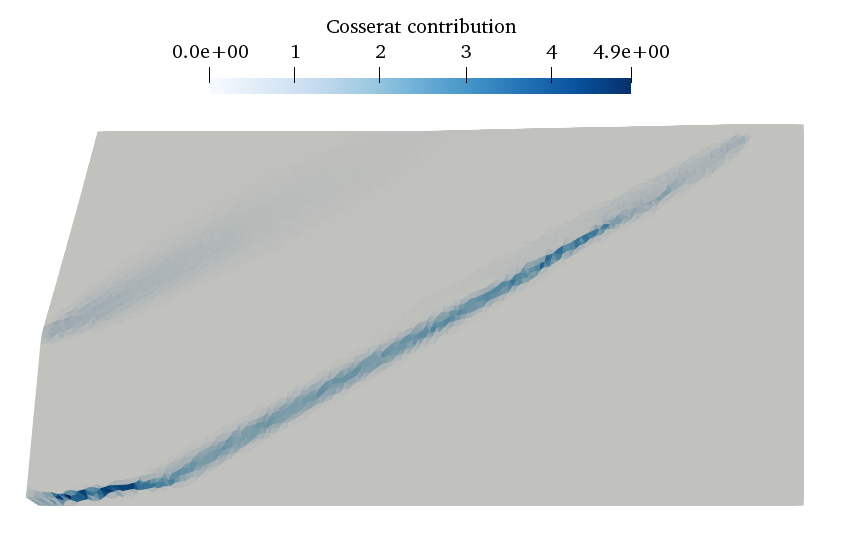}
\caption{$\ell/H=0.01$}
\end{subfigure}
\hfill
\begin{subfigure}{0.32\textwidth}
\includegraphics[width=\textwidth]{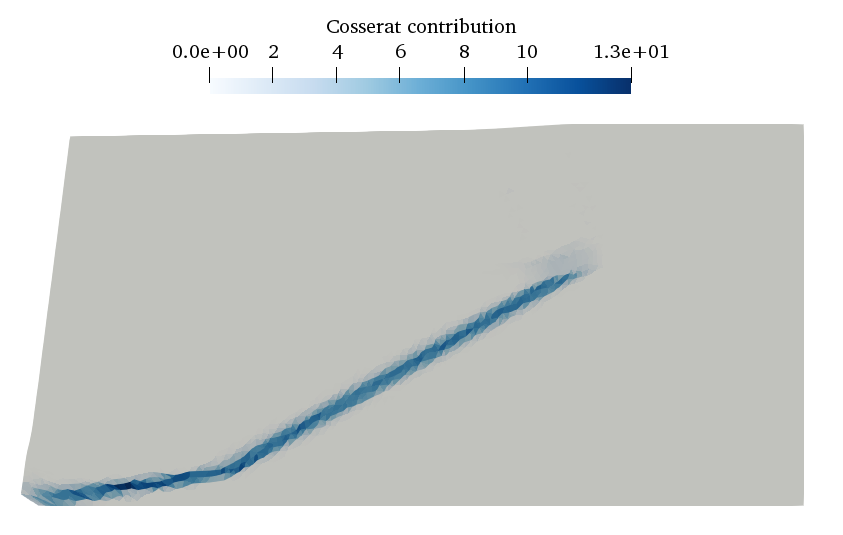}
\caption{$\ell/H=0.1$}
\end{subfigure}
\end{center}
\caption{Collapse mechanism of the jointed rock excavation and pure Cosserat contribution to the total dissipation}
\label{cosserat-mechanisms}
\end{figure}

A mixed approach for this model has been implemented in the spirit of \eqref{mixed-LA} which avoids the need to compute the support function expression associated with \eqref{joint-rock-G}. Continuous $\PP^2$ (resp. $\PP^1$) Lagrange elements have been used for $\bu$ (resp. $\omega$) and discontinuous $\PP^1_d$-Lagrange elements for both $\bSig$ and $\boldsymbol{H}$. The strength conditions have been imposed at the vertices of each element. We considered the stability problem of an excavation of height $H$, making a $25^{\circ}$ angle with the vertical and subjected to its self-weight of intensity $\gamma$. The problem amounts to find the maximum value of the non-dimensional stability factor $K^+=\left(\dfrac{\gamma H}{c_m}\right)^+$. For numerical applications, we took $c_j=0.5c_m$, $\phi_j=20^{\circ}$, $\phi_m=40^{\circ}$, $\theta=10^{\circ}$ and varied the joint spacing $\ell$. The evolution of the stability factor estimates as a function of $\ell/H$ has been represented in Figure \ref{cosserat-results} for two different mesh sizes. As for the strain gradient model, strengthening is observed for increasing $\ell/H$ ratios. Interestingly, size-effects are much stronger for this problem than those of Figure \ref{strain-gradient-results}. The obtained value in the standard Cauchy ($\ell=0$) case is quite close to the analytical upper bound of $K^+\leq 1.47$ derived for the same problem in \cite{de2002failure}. Finally, collapse mechanisms and a measure of the pure Cosserat contribution $(\Sigma_{21}-\Sigma_{12})(u_{2,1}-u_{1,2} - \omega)+\boldsymbol{H}\cdot\nabla \omega$ to the total dissipation have been represented in Figure \ref{cosserat-mechanisms}. It can be observed that the shape of the collapse mechanism and the location of "shearing" zones involving Cosserat effects is quite dependent on the joint spacing. For $\ell=0$, a triangular sliding block with a concentrated slip zone is obtained, approximately corresponding to the merging of the two slip bands of Figure \ref{cosserat-mech-1e-3}.

\section{Conclusions}

This paper proposed a way to easily formulate and solve limit analysis problems by taking advantage of three essential ingredients:
\begin{itemize}
 \item the FEniCS finite element library with its high-level domain specific language and large choice of finite-element interpolations
 \item the representation of limit analysis criteria and associated support functions in a conic programming format
 \item the resolution of the corresponding conic optimization problems by a dedicated and efficient interior-point solver implemented in the \texttt{Mosek} software package
 \end{itemize} 
The first two ingredients offer an extremely large versatility in the problem definition, giving access to an extremely large range of applications. In particular, various finite-element interpolations can be easily defined, giving access to traditional continuous displacement-based upper-bound formulations but also their discontinuous counterpart through Discontinuous Galerkin function spaces. Equilibrium lower-bound elements are therefore also accessible through such spaces as well as mixed formulations through different choices of quadrature rules. We also showed that plates and shells problems could also be discretized without effort. The conic representation format \eqref{conic-representable} offers a unified way of defining strength conditions and associated support functions for different mechanical models, ranging from classical 2D/3D continuum mechanics to plate bending criteria including potential shear conditions, shell criteria with membrane/bending interaction or even generalized continua such as Cosserat or strain gradient models. If the conic representation format is large enough to encompass many strength criteria, it is also sufficiently disciplined to yield optimization problems of the conic programming class for which dedicated solvers like \texttt{Mosek} have been designed. \texttt{Mosek} is indeed known to be a state-of-the-art optimizer for this class of problems and therefore offers efficiency and robustness of the solution procedure.

Obviously, the present work could still improve upon some aspects, in particular regarding computational efficiency. For instance, many additional auxiliary variables are usually introduced when complying with the conic programming format. Some of them may be handled and eliminated by \texttt{Mosek} during its pre-processing phase, although this is not entirely clear since it is used as a black-box. Devising an interior-point solver specific to limit analysis problems can take advantage of the problem structure and may be more efficient. 
Besides, there are two points which prevent solving extremely large-scale 3D problems. The first one is related to the use of direct solvers in the interior-point inner iterations which requires large memory capacities for 3D problems. The use of iterative solvers is still an active research topic due to the difficulty of efficiently preconditioning the interior-point linear systems. The second one concerns the need to solve SDP problems when SDP-representable criteria like Mohr-Coulomb, Tresca or Rankine are used in 3D. Even though interior-point solvers efficiency has greatly improved for SDP problems over the last decade, it is still more difficult to solve than an SOCP problem of similar size. Improving even more their efficiency or finding alternate strategies would be a great benefit for such 3D problems.

If mathematical programming (and more particularly conic programming) tools for solving limit analysis has now emerged as the state-of-the art method, some extensions have already been proposed in the literature to apply them also to closely related problems. One can mention, for instance, elastoplasticity \cite{krabbenhoft2007interior,krabbenhoft2007formulation}, viscoplasticity for yield stress fluids \cite{bleyer2015efficient,bleyer2018advances}, contact in granular materials \cite{krabbenhoft2012granular,zhang2014particle}, limit analysis-based topology optimization \cite{kammoun2014direct, herfelt2019strength}, etc. The present framework is sufficiently general to also extend to these related problems, see for instance the application to viscoplastic fluids in \cite{bleyer2019automating}. Other situations appear however more difficult to include such as non-associative behaviours or geometrical non-linearities since such problems cannot be formulated as convex optimization problems anymore. Nonetheless, some works have already proposed some iterative strategies for tackling non-associativity \cite{gilbert2006limit,krabbenhoft2012associated,portioli2014limit}, it would therefore be interesting to pursue in this direction.

\bibliographystyle{elsarticle-harv}
\bibliography{fenics_limit_analysis}
\end{document}